\documentclass[hidelinks,11pt]{siamart171218}

\usepackage[a4paper, margin=1.1in]{geometry}

\usepackage{amssymb, amsmath}
\usepackage{mathtools}
\usepackage{lineno}
\usepackage[font=small,labelfont=bf]{caption}
\usepackage{bbm}
\usepackage{tikz}
\usepackage{float}

\makeatletter
\newcommand\RedeclareMathOperator{%
  \@ifstar{\def\rmo@s{m}\rmo@redeclare}{\def\rmo@s{o}\rmo@redeclare}%
}
\newcommand\rmo@redeclare[2]{%
  \begingroup \escapechar\m@ne\xdef\@gtempa{{\string#1}}\endgroup
  \expandafter\@ifundefined\@gtempa
     {\@latex@error{\noexpand#1undefined}\@ehc}%
     \relax
  \expandafter\rmo@declmathop\rmo@s{#1}{#2}}
\newcommand\rmo@declmathop[3]{%
  \DeclareRobustCommand{#2}{\qopname\newmcodes@#1{#3}}%
}
\@onlypreamble\RedeclareMathOperator
\makeatother
\RedeclareMathOperator{\div}{div}

\DeclareMathOperator{\bfx}{\boldsymbol{x}}

\DeclareMathOperator{\kbfx}{\widehat{\boldsymbol{x}}}
\DeclareMathOperator{\td}{d}

\DeclareMathOperator{\dt}{d{\it t}}
\DeclareMathOperator{\ds}{d{\it s}}
\DeclareMathOperator{\knabla}{\widehat{\nabla}}

\DeclareMathOperator{\dbfx}{d\boldsymbol{x}}
\DeclareMathOperator{\dkbfx}{d\widehat{\boldsymbol{x}}}
\DeclareMathOperator{\dS}{dS}

\DeclareMathOperator{\kx}{\widehat{\it x}}
\DeclareMathOperator{\ky}{\widehat{\it y}}

\DeclareMathOperator{\kdiv}{\widehat{div}}

\DeclareMathOperator{\bfw}{\boldsymbol{w}}

\DeclareMathOperator{\bfn}{\boldsymbol{n}}

\DeclareMathOperator{\bfomega}{\boldsymbol{\omega}}
\DeclareMathOperator{\bfpsi}{\boldsymbol{\psi}}
\DeclareMathOperator{\kbfpsi}{\widehat{\boldsymbol{\psi}}}
\DeclareMathOperator{\kbfomega}{\widehat{\boldsymbol{\omega}}}

\DeclareMathOperator{\kbfw}{\widehat{\boldsymbol{w}}}

\DeclareMathOperator{\kbfu}{\widehat{\boldsymbol{u}}}
\DeclareMathOperator{\kbfg}{\widehat{\boldsymbol{g}}}

\DeclareMathOperator{\kf}{\widehat{\it f}}
\DeclareMathOperator{\kg}{\widehat{\it g}}
\DeclareMathOperator{\ku}{\widehat{\it u}}

\DeclareMathOperator{\kpsi}{\widehat{\psi}}
\DeclareMathOperator{\jac}{\cal J}

\DeclareMathOperator{\kALE}{{\widehat{\cal A}}}
\DeclareMathOperator{\kALEF}{{\widehat{\cal F}}}

\DeclareMathOperator{\kALEinvt}{{\widehat{\cal A}^{\textit{-1}}_{\it t}}}
\DeclareMathOperator{\kjac}{\widehat{\cal J}}

\DeclareMathOperator{\kF}{\widehat{\it F}}
\DeclareMathOperator{\tA}{\mathbbm{A}}
\DeclareMathOperator{\tB}{\mathbbm{B}}

\DeclareMathOperator{\mR}{\mathbb{R}}

\DeclareMathOperator{\kOmega}{\widehat{\Omega}}
\DeclareMathOperator{\kQ}{\widehat{\it Q}}

\DeclareMathOperator{\kK}{\widehat{\it K}}

\newsiamthm{remark}{Remark}

\begin{document}

\title{ALE--type FEM formulation for PDE\MakeLowercase{s} on time--dependent domains with vanishing discrete SCL\thanks{Submitted to the editors \today .\funding{This work was supported in part by NHRI (National Health Research Institutes, Institute of Biomedical Engineering and Nanomedicine, Zhunan, Taiwan, project number BN-106-PP-08) and MOST (Ministry of Science and Technology, Taipei, Taiwan, project number MOST--106--2115--M--400--001--) project.}}}


\author{
	Filip Ivan\v ci\' c\footnotemark[2] \footnotemark[3]  
	\and Tony W. H. Sheu \footnotemark[2] \footnotemark[4] \footnotemark[5] \footnotemark[6] 
	\and Maxim Solovchuk\footnotemark[2] \footnotemark[3] \footnotemark[6] 
}
\maketitle

\renewcommand{\thefootnote}{\fnsymbol{footnote}}

\footnotetext[2]{Department of Engineering Science and Ocean Engineering, National Taiwan University, Taipei, Taiwan}
\footnotetext[3]{National Health Research Institutes, Institute of Biomedical Engineering and Nanomedicine, Zhunan, Taiwan}
\footnotetext[4]{Center for Advanced Study in Theoretical Sciences (CASTS), National Taiwan University, Taipei, Taiwan} 
\footnotetext[5]{Institute of Applied Mathematical Sciences, National Taiwan University, Taipei,
Taiwan}
\footnotetext[6]{Corresponding author.}
\footnotetext[0]{{\it email addresses:} ivancicfi@gmail.com (Filip Ivan\v ci\' c), twhsheu@ntu.edu.tw (Tony W. H. Sheu), solovchuk@gmail.com (Maxim Solovchuk)}
\renewcommand{\thefootnote}{\arabic{footnote}}

\begin{abstract}
The aim of this paper is to introduce a finite element formulation within Arbitrary Lagrangian Eulerian framework with vanishing discrete {\it Space Conservation Law} (SCL) for differential equations on time dependent domains. The novelty of the formulation is the method for temporal integration which results in preserving the SCL property and retaining the higher order accuracy at the same time. Once the time derivative is discretized (based on integration or differentiation formula), the common approach for terms in differential equation which do not involve temporal derivative is classified to be a kind of "time averaging" between time steps. In the spirit of classical approaches, this involves evaluating these terms in several points in time between the current and the previous time step ($[t_n,t_{n+1}]$), and then averaging them in order to provide the satisfaction of discrete SCL. 
Here, we fully use the polynomial in time form of mapping through which evolution of domain is realized -- the so called ALE map -- in order to avoid the problematics arising due to the moving grids. We give a general recipe on temporal schemes that have to be employed once the discretization for the temporal derivative is chosen. Numerical investigations on stability, accuracy and convergence are performed and the simulated results are compared with benchmark problems set up by other authors.
\end{abstract}



\section*{Introduction}\label{sec:introduction}
When dealing with problems in time--dependent domains, usually in the context of fluid dynamics, the Arbitrary Lagrangian Eulerian (ALE) framework is probably the most natural set-up environment (certainly among the most popular ones) in finite element and finite volume methods. A detailed description of the ALE framework can be found in \cite{donea} and references therein. ALE formulation is realized through a homeomorphic map (usually referred to as the {\it ALE map}), which at each time $t$ associates current configuration (a physical domain) to some fixed referent one. This ALE map is prescribed (or can somehow be obtained for current time) on the domain boundary, while it is quite arbitrary in interior. Consequently, there is some freedom in defining it in the domain interior. Several issues have been raised by various authors on the regularity of ALE map in order to ensure the problem of interest is well posed, as well as on the concrete discrete realizations of it. It has been noticed that some properties of the ALE map, which are trivially satisfied on the continuous level, fail on the discrete level and introduce, in turn, extra errors in numerical approximation of solution as well as influence the level of stability of the resulting algorithm. This issue was noticed with somewhat delay compared to that with the introduction of ALE formulations. Moreover, many authors do not pay any special attention to it. However, in recent years, numerous studies have been published reporting conspicious results and issues on this topic  in the context of finite element, finite volume and finite difference methods, and the core of them is the so called {\it Space Conservation Law} ({\it SCL}) introduced in \cite{trulio,thomas-lombard}. The heart of this issue (which will be formally defined later) is the method for the temporal integration of the considered equations and the effects arising from the indispensable discrete grid velocity. Most of the original work on this topic has been done in context of the {\it finite difference method} (\cite{thomas-lombard, abe,visbal,deng,yee}) and {\it finite volume method} (\cite{zhang,demirdzic,farhat96,farhat99,farhat00,farhat01}), and, lately, extended to the {\it finite element method} (\cite{farhat96,farhat99,formaggia-nobile99,formaggia-nobile04,boffi-gastaldi,etienne}). In this work, finite element method is considered.

In the context of fluid dynamics and related problems, these mentioned effects usually appear in the form of artificial sinks and sources and may possibly influence stability and accuracy of the numerical method. The ALE formulation consists of rewriting the equations, naturally posed in the moving domain, with respect to fixed referential configuration, employing the {\it Reynolds transport theorem} in the process. In this framework, the temporal derivative expressed with respect to referent configuration introduces an extra term in which grid velocity appears and is known as the cause leading to the mentioned issues. 

The attempt of this paper is to intertwine the topics from above, ALE framework and SCL, from, somehow, different perspective than majority of papers we came across up to this point. The problematics with grid velocity in ALE approach for treating the moving domain problems appears in the transition from continuous to discretized formulation (in time). The mentioned SCL condition is a trivial identity on continuous level which only comes to the fore in discrete counterpart of the considered problem. In its integral version it takes the form (\cite{thomas-lombard})\[ \frac{\td}{\dt}\int_{\Omega(t)}\dbfx=\int_{\Omega(t)}\div\bfw\dbfx, \] which is also often referred to as the {\it finite volume form}, with $\Omega(t)$ being the domain at time \(t\) and $\bfw$ the domain velocity. There are several different versions of the {\it finite element form} of SCL, but essentially they are all derived from the above identity. It seems that the usual approach in dealing with mentioned numerical issues is based on employing a better time integrator. For example, Formaggia and Nobile in \cite{formaggia-nobile99,formaggia-nobile04} proposed multi-point time quadrature of the right--hand side of the system of {\it ordinary differential equations} (ODE's) in time chosen in a way such that the SCL is not violated: to solve ODE of form $dy/dt=f(t,y(t))$, they replace the usual backward Euler  method \[ y^{n+1}-y^n=\Delta t f(t_{n+1},y^{n+1}) \textrm{\space\space\space by \space\space\space}y^{n+1}-y^n=\Delta t \sum_{l=0}^m w_l f(t^l_{n,n+1},y^{n+1}),\] $t_{n,n+1}^l$ being the time instant in $[t_n,t_{n+1}]$ and $w_l$ appropriately chosen {\it weights}. For large finite element problems this approach seems quite an expensive one from computational standards, especially if all quantities are evaluated in the current configuration. In that case for each involved time instant, finite element space has to be updated due to its dependence on time through the domain.

In \cite{boffi-gastaldi} a similar approach as in \cite{formaggia-nobile04} is taken. They investigate a model ALE scheme with respect to different choices of time discretizations. For considered scheme they investigate the relationships between the scheme stability and the SCL. In both \cite{formaggia-nobile04,boffi-gastaldi} they show that satisfying the SCL condition is neither a necessary, nor a sufficient condition for stability, except for the implicit Euler scheme. The same conclusion is derived in this paper as well.

Here, however, the full potential of the polynomial (in time) form of ALE map will be used. Since the polynomials can be integrated exactly, by ably transforming the problem into some appropriate configuration, the SCL can be made trivially satisfied again in the discrete counterpart of problem formulation. Apart from its very natural and intuitive look, we will show that adapting this approach to a chosen temporal discretization is straightforward. Generalizations to the higher order schemes based on some integration or differentiation formulas are straightforward to derive. Also, there is more room in choosing the way of calculation of the grid velocity without (essentially) changing the formulation of the problem, thus making the coding easily adjustable. Implementationwise, everything is kept on the original referent configuration and consequently introduces some additional differential operators in order to transform the space derivatives from current configuration onto the referent one. We believe this approach offers a few advantages ahead of the implementation on current configuration (of course, under the assumption that the two approaches are equivalent on the discrete level which will be argued later). Among them, we mention three here. The first two lie in the fact that the referential domain is fixed in time. Then the test/basis function spaces will be time independent and the finite element spaces do not need to be updated at each time step -- a huge advantage from computational point of view. The third, keeping everything on the referential configuration, the evolution of the domain is kept in the Jacobians of the ALE map thus making the connection between the domain time--dependency and all the terms in differential equation more clear. In cases when a weak formulation is posed on the current configuration, the only explicit connection between the domain time--dependency and under integral terms is through domain velocity which appears in the form of convective terms.

The paper is organized as follows. In Section~\ref{sec:ALE} the ALE formulation is reviewed and the notation is introduced. In Section~\ref{sec:cootransform} the pullback to referent configuration which is more convenient for working with problems on time--dependent domains is illustrated on heat equation. Section~\ref{sec:SCL} deals with space conservation law -- the basis for the introduction of the new type formulation. In Section~\ref{sec:vanishSCL}, the new formulation with vanishing discrete SCL is introduced and two different ways of calculating the grid velocity are presented. In Section~\ref{sec:schemes} some of the most popular schemes (implicit Euler and Crank--Nicolson schemes, backward differentiation formulas) are adapted in the environment of our new formulation. In Section~\ref{sec:fem}, some insights on the restriction of the spatial discretization of the ALE map with finite elements are given. Section~\ref{sec:numerics} deals with numerical validation. Finally, in Section~\ref{sec:conclusion} we give the summary of what has been done and draw the conclusions.

\section{The Arbitrary Lagrangian Eulerian formulation}\label{sec:ALE}
In this section, the arbitrary Lagrangian Eulerian (ALE) formulation for moving domain problems is reviewed. All the details on the results here and more mathematically rigorous derivations of them can be found in \cite{donea,formaggia-nobile99,formaggia-nobile04} and references therein. 

The general idea of the ALE framework consists of the interplay between the fixed referential domain (which most often coincides with initial domain, but does not have to) and the current physical domain occupied by the medium. The interplay between these two domains is realized through the so called {\it ALE map} which maps the referential domain into the current one. In order to perform the necessary calculus, a minimal smoothness of the ALE map  has to be demanded  -- e.g. a kind of inverse of the ALE map has to exist in order to ensure the correspondence between the referential domain and the current one (in sense that one can be obtained from the other). 

Let $\kOmega\subset\mR^d$, $d=2,3,$ be a fixed referential domain and $\Omega(t)\equiv\Omega_t\subset\mR^d$ the current (physical) domain occupied by the fluid. It is assumed that boundaries of the domain are sufficiently smooth -- this usually refers to the Lipschitz continuous boundary -- and that the domain evolution can be followed through a one--parameter family of mappings $(\kALE_t)_{t\in[0,T]}\equiv(\kALE(\cdot,t))_{t\in[0,T]}$, $T<\infty$,
\begin{equation}
\begin{split}
\kALE_t\colon\kOmega\to \mR^d\textrm{, }t\in[0,T] \textrm{\space\space , \space\space}(\kbfx,t)\mapsto(\bfx,t)\textrm{, }\kbfx\in\kOmega\textrm{, }\bfx\in \Omega(t),
\end{split}
\end{equation}
i.e. $\kALE_t$ maps the referential into the current domain, $\kOmega\mapsto\Omega_t\equiv\Omega(t)=\kALE_t(\kOmega)$. In this context, we refer to $\kbfx\in\kOmega$ as the {\it ALE coordinate} while to $\bfx=\kALE_t(\kbfx)\in\Omega(t)$ we refer as an {\it Eulerian (or spatial) coordinate}.

Denote by $\kQ_T=\kOmega\times(0,T)$ the referential domain and by \[Q_T=\bigcup_{t\in(0,T)}\Omega(t)\times\{ t\}\] the domain trajectory (sometimes we write $Q_T=\Omega(t)\times(0,T)$). Let $f\colon Q_T\to\mR$ and $\kg\colon\kQ_T\to\mR$ be two scalar fields defined on the current and  the referential configurations, respectively. We define their {\it ALE} and {\it Eulerian} counterparts respectively by
\begin{equation}
\begin{split}
\kf\colon\kQ_T\to\mR\textrm{, }\kf=f\circ\kALE_t \textrm{ \space\space and \space\space }g\colon Q_T\to\mR\textrm{, } g = \kg\circ\kALEinvt.
\end{split}
\end{equation}

In the ALE environment, the temporal derivative of an Eulerian field can be considered from different viewpoints. Let $f\colon Q_T\to\mR$ be an Eulerian field, and $\kf=f\circ\kALE_t$  its ALE counterpart. The time derivative of an Eulerian field $f$ in the ALE frame -- i.e. time derivative in the ALE frame, written with respect to the spatial coordinates -- is defined as
\begin{equation}\label{ALEtimederivative}
\frac{\partial}{\partial t}\Big|_{\kbfx}f\colon Q_T\to\mR\textrm{ , }\frac{\partial }{\partial t}\Big|_{\kbfx}f(\bfx,t) = \frac{\partial\kf}{\partial t}(\kbfx, t)\textrm{, }\kbfx=\kALEinvt(\bfx).
\end{equation}
The time derivative of an Eulerian field in the spatial frame is just the classical time partial derivative \[\frac{\partial }{\partial t}\Big|_{\bfx}f = \frac{\partial f}{\partial t}.\] 
Now, the domain velocity is defined as
\begin{equation}
\bfw(\bfx,t) = \frac{\partial}{\partial t}\Big|_{\kbfx}\bfx\textrm{, }\bfx=\kALE(\kbfx,t)\textrm{ , i.e. } \bfw(\bfx,t) = \kbfw(\kbfx,t)=\frac{\partial}{\partial t}\kALE(\kbfx,t).
\end{equation}

In terms of displacement, ALE map can be written in the form of
\begin{equation}
\kALE(\kbfx,t)=\kbfx+\kbfu(\kbfx,t) \textrm{\space\space , \space\space}\kbfu(\kbfx,t)=\bfx(\kbfx,t)-\kbfx,
\end{equation}
$\kbfu(\kbfx,t)$ being the displacement of $\kbfx$ at time $t$.
By application of the chain rule it is straightforward to obtain
\begin{equation}\label{alederivative}
\begin{split}
\frac{\partial}{\partial t}\Big|_{\kbfx}f & = \frac{\partial}{\partial t}\Big|_{\bfx}f + \frac{\partial \bfx}{\partial t}\Big|_{\kbfx}\cdot\nabla f = \frac{\partial f}{\partial t} + \bfw\cdot\nabla f.
\end{split}
\end{equation}
As a result, by introducing an ALE temporal derivative instead of an Eulerian one, an extra convective--type term, due to domain movement, will be introduced in the equation of interest.

An important role in interplay between configurations play the gradient and the Jacobian of ALE map, given respectively by 
\begin{equation}
\kALEF_t=\knabla\kALE_t\textrm{ , }\kjac_t = \det(\knabla\kALE_t).
\end{equation} 
The so called {\it Euler expansion formula} states
\begin{equation}\label{eulerexpansion}
\frac{\partial}{\partial t}\kjac_t = \kjac_t\widehat{\div\bfw}.
\end{equation}

\subsection{Test function spaces in ALE framework}\label{sec-test-func}
Finite element method is based on weak formulation of the considered partial differential equation in which space of test functions plays an essential role. The standard approach on stationary domains is to take test functions independent of time, and the similar approach is followed in current context. However, since here we are dealing with time--dependent domain, following this arrangement, a whole family of test function spaces is obtained such that for each $t$ test functions on $\Omega_t$ are time independent (in the context of the current configuration $\Omega_t$).  Surely, this deserves some special attention.

Let $V(\kOmega)$ be a space of admissible test functions defined on a referent domain which consists of regular enough functions $\kpsi\colon\kOmega\to\mR$. Admissible here stands for {\it well defined} in the sense of preserving boundary conditions -- e.g., $\kpsi\in V(\kOmega)$ has to vanish on Dirichlet part of the boundary so $V(\kOmega)$ depends on the problem itself. By regular enough we mean that all integrals in weak formulation of the problem make sense -- e.g. after partial integration and transferring derivatives onto test function, all terms under integral sign are integrable. Usually, for second order PDEs, $V(\kOmega)$ will be some subset of $H^1(\kOmega)$, with $H^1$ being the {\it Sobolev space} here defined in the standard way (\cite{brezis}). Then, through ALE mapping, we can identify a corresponding set $V(\Omega_t)$ of admisible test functions on the current configuration defined as
\begin{equation}
V(\Omega_t) = \{ \psi\colon\Omega_t\times[0,T)\to\mR\mid \psi=\kpsi\circ\kALEinvt\textrm{, }\kpsi\in V(\kOmega) \}.
\end{equation}
Taking into account the correspondence between $V(\kOmega)$ and $V(\Omega_t)$ and the definition of ALE temporal derivative $\displaystyle\frac{\partial}{\partial t}\Big|_{\kbfx}$, the following relation is obtained
\begin{equation}
0 = \frac{\partial\kpsi}{\partial t} = \frac{\partial}{\partial t}\Big|_{\kbfx} \psi = \frac{\partial \psi}{\partial t} + \bfw\cdot\nabla\psi\textrm{, }\forall \psi\in V(\Omega_t).
\end{equation} 
Now, let $f=f(\bfx,t)$ be an arbitrary time--differentiable Eulerian field. By employing the chain rule, it follows
\begin{equation}
\frac{\partial(\psi f)}{\partial t}\Big|_{\kbfx} = \psi\frac{\partial f}{\partial t}\Big|_{\kbfx}\textrm{, }\forall \psi\in V(\Omega_t).
\end{equation}
Recalling the Reynolds transport theorem, the following identities are obtained
\begin{equation}\label{reyn1}
\frac{\td}{\dt}\int_{\Omega_t}\psi\dbfx = \int_{\Omega_t}\psi\div\bfw\dbfx,
\end{equation}
\begin{equation}\label{propertydtintpsiu}
\frac{\td}{\dt}\int_{\Omega_t}\psi f = \int_{\Omega_t}\psi\left( \frac{\partial  }{\partial t}\Big|_{\kbfx}f + f\div\bfw \right)\dbfx,
\end{equation}
\begin{equation}\label{continuousgcl}
\frac{\td}{\dt}\int_{\Omega_t}\psi \chi = \int_{\Omega_t}\psi\chi\div\bfw,
\end{equation}
for all $\psi,\chi\in V(\Omega_t)$.

\subsection{A construction of ALE map}
In practical problems, one usually faces the problem of determining the ALE map $\kALE_t$ in $\kOmega\times(0,T)$ with $\kALE_t$ being prescribed on the boundary $\partial\kOmega\times(0,T)$, i.e.
\begin{equation}
\begin{split}
&\textrm{given }\kALE_t(\kbfx)=\kbfg(\kbfx,t)\textrm{ on }\partial\kOmega\times(0,T)\\
&\textrm{determine }\kALE_t\textrm{ in }\kOmega\times(0,T).
\end{split}
\end{equation}
Several techniques have been proposed in the literature for the extension of $\kALE_t$ from the boundary into the interior. Based on the problem we are facing, one technique may be more appropriate than the other. Classically, one solves a PDE for unknown $\kALE_t$ in $\kOmega\times(0,T)$ subjected to Dirichlet boundary conditions. Since most often ALE map is needed only at discrete time levels, this adds up to solving a PDE for unknown $\kALE_{n+1}$ with  $\kALE_{n+1}=\kbfg_{n+1}$ prescribed on $\partial\kOmega$, where $\kALE_{n+1}=\kALE(\cdot,t_{n+1})$. Once $\kALE_{n+1}$ is known, $\Omega_{n+1}=\Omega(t_{n+1})$ can be obtained. The most popular extension in the literature seems to be a harmonic extension but we mention that alternative approaches are possible. For example, instead of the Laplace equation (harmonic extension), the linear elasticity equation is a quite popular choice as well.

\section{The transformation of configurations}\label{sec:cootransform}
The SCL non--violating formulation we propose here turns out to be particularly convenient if the PDE of interest is pulled back to the referent configuration. However, during the {\it pullback} procedure, spatial differential operators (originally functions of current configuration) have to be transformed into their respective counterparts in the referent configuration. In this section we summarize the necessary results, illustrate the pullback procedure on the simple heat equation, and recall some "tricks" which make the transition from the current to the referent configuration profitable.

\subsection{Pullback of the heat equation}
Consider a heat equation subjected to the initial and Dirchlet boundary conditions
\begin{equation}\label{simpleheat}
\begin{split}
\frac{\partial}{\partial t}u - \alpha\Delta u & = f \textrm{ \dots\space }\Omega(t)\times(0,T),\\
u(0) & = u_0\textrm{ \dots\space }\Omega(0),\\
u & = u_D\textrm{ \dots\space }\partial\Omega(t)\times(0,T),
\end{split}
\end{equation}
with $\alpha$ being a constant for simplicity and assuming $(\kOmega,t)\mapsto\Omega(t)$ is realized through $\kALE _t$ which is prescribed. Employing the ALE time derivative (\ref{simpleheat})$_1$ becomes
\begin{equation}\label{strong-heat}
\frac{\partial}{\partial t}\Big|_{\kbfx} u - \alpha\Delta u -\bfw\cdot\nabla u = f \textrm{ \dots\space }\Omega(t)\times(0,T),
\end{equation}
and the corresponding conservative weak formulation is found
\begin{equation}\label{simpleheatweak}
\frac{\td}{\dt}\int_{\Omega(t)}\psi u \dbfx +\int_{\Omega(t)}\left[ \alpha\nabla\psi\cdot\nabla u - \psi\bfw\cdot\nabla u - \psi u \div\bfw -\psi f\right]\dbfx = 0
\end{equation}
with an appropriate choice of the test function $\psi$ (e.g. $\psi$ has to vanish on the $\partial\Omega(t)\times(0,T)$ due to Dirichlet boundary conditions). Transforming the weak formulation (\ref{simpleheatweak}) onto the referent configuration, (\ref{simpleheatweak}) takes the following form that we will work with in this paper:
\begin{equation}\label{heatweakreferent}
\begin{split}
0 =  \frac{\td}{\dt}\int_{\kOmega} \kpsi\ku & \kjac_t\dkbfx \\
& + \int_{\kOmega}\left[ \alpha\widehat{\nabla\psi}\cdot\widehat{\nabla u}\kjac_t - \kpsi\kbfw\cdot\widehat{\nabla u}\kjac_t - \kpsi\ku\widehat{\div\bfw}\kjac_t - \kpsi\kf\kjac_t \right]\dkbfx.
\end{split}
\end{equation}
\begin{remark}
One should pay an extra attention on terms $\widehat{\nabla\psi}$, $\widehat{\nabla u}$, $\widehat{\div\bfw}$ and notice that the differential operators in these terms operate with respect to the current configuration and then the whole expression is pulled back to the referent configuration, e.g.
\[ \widehat{\div\bfw} = (\div\bfw)\circ\kALE_t\neq \kdiv\kbfw. \] This results in a problem, since the current configuration is not necessarily known and might be part of the problem as well. Luckily, thanks  to the ALE map, all of these operators can be expressed with respect to the referent configuration.
\end{remark}

\subsection{Transformation of first order spatial differential operators}
We wish to express the differential operators defined in physical configuration into their referent configuration counterparts. To begin with, note that by the chain rule
\begin{equation}
\frac{\partial}{\partial x_i} = \frac{\partial\kx_j}{\partial x_i}\frac{\partial}{\partial\kx_j},
\end{equation}
with  {\it Einstein summation convention}\footnote{In Einstein summation convention there is an implied summation over the terms with the repeated index.} being employed for the sake of compact presentation. Now, by the direct application of above chain rule, it is easy to derive
\begin{equation}
\nabla\psi = \nabla\kpsi = \kALEF_t^{-T}\knabla\kpsi \textrm{ and }\div\bfpsi = \div\kbfpsi = \frac{1}{\kjac_t}\kdiv(\kjac_t\kALEF_t^{-1}\kbfpsi)
\end{equation}
for scalar field $\psi$ and vector field $\kbfpsi$.

Note that in the above expressions $\kALEF_t^{-1}$ appears which is a function of $\bfx$ i.e. defined on physical configuration, while our framework is set up on the referent configuration. To express $\kALEF_t^{-1}$ in the referent framework, we recall some facts from linear algebra. Assuming that we are given regular matrices $\tA=(a_{ij})$ and $\tB=(b_{ij})$ such that $\tB=\tA^{-1}$, it is well known that 
\begin{equation}
b_{ij} = \frac{1}{\det\tA}(-1)^{i+i}M_{ji}
\end{equation}
where $(-1)^{i+j}M_{ij}$ is the $(i,j)$--cofactor of $\tA$. 

In the spirit of the above discussion, we can now transform $\kALEF_t^{-1}$ into its counterpart defined on the referent configuration. For the sake of illustration we show this in two dimensions, while the three dimensional case is similar. Using the notation $\kbfx=\kALE_t^{-1}(\bfx)$ (due to invertibility of $\kALE_t$ ) we have
\begin{equation}
\kALEF_t^{-1} = \begin{bmatrix}
\partial_x\kx & \partial_y\kx \\
\partial_x\ky & \partial_y\ky
\end{bmatrix} = \frac{1}{\kjac_t}\begin{bmatrix}
\partial_{\ky}y & -\partial_{\ky}x \\
-\partial_{\kx}y & \partial_{\kx}x
\end{bmatrix} = \frac{1}{\kjac_t}\kF_t
\end{equation}
where
\begin{equation}
\kF_t \coloneqq  \begin{bmatrix}
\partial_{\ky}y & -\partial_{\ky}x \\
-\partial_{\kx}y & \partial_{\kx}x
\end{bmatrix}.
\end{equation}

Finally, we have the necessary identities for scalar and vector fields $\psi$ and $\bfpsi$ being transformed onto the referent configuration:
\begin{equation}
\begin{split}
\widehat{\nabla\psi}  = \frac{1}{\kjac_t}\kF_t^T\knabla\kpsi \textrm{ and }\widehat{\div\bfpsi}  = \frac{1}{\kjac_t}\kdiv(\kF_t\kbfpsi).
\end{split}
\end{equation}
Employing them enables us to rewrite (\ref{heatweakreferent}) into the form convenient for the numerical implementation:
\begin{equation}
\begin{split}
0=\frac{\td}{\dt} & \int_{\kOmega}  \kpsi\ku\kjac_t  \dkbfx \\
+ & \int_{\kOmega}\left[ \alpha\frac{1}{\kjac_t}\kF_t\kF_t^T\knabla\kpsi\cdot\knabla\ku - \kpsi\kF_t\kbfw\cdot\knabla\ku  - \kpsi\ku\kdiv(\kF_t\kbfw) - \kpsi\kf\kjac_t\right]\dkbfx
\end{split}
\end{equation}
\begin{remark} In context of the spatial differential operators a question whether the discrete derivatives on the current and their counterparts on the referential domain coincide may arise. However, this is not expected to be a problem if the ALE map is polynomial in space (since then the grid is moved by polynomial) and polynomials are treated exactly.
\end{remark}
\begin{remark} From the computational point of view, it seems that pulling back to the referential configuration may offer some advantages. In problems posed on moving domains (fluid structure interaction, free surface problems, etc.) the movement of the domain is itself part of the problem. Very common solution involves the algorithms of iterative nature based on the fixed point theorems. In context of finite elements, this means that at each iteration the ALE map is updated, grid adjusted accordingly and finite element spaces have to be updated --  computationally a quite expensive procedure. Keeping the environment on the referential configuration only the ALE map has to be updated while the finite element spaces remain unchanged. 
\end{remark}
\section{The Space Conservation Law}\label{sec:SCL}
As far as authors know, the first mentioning of the {\it Space conservation law} (SCL) was by Trulio and Trigger in \cite{trulio}. SCL was latter rediscovered and formalized by Thomas and Lombard in \cite{thomas-lombard}, and its differential form was derived. A significant work on the subject has been done in the series of papers \cite{farhat96,farhat99,farhat00,farhat01}. In \cite{farhat96} they showed that for {\it space--time finite elements} (STFEM) the discrete SCL vanishes. However, since for STFEM (\cite{hughes-stfem}) the dimension of the discrete system to solve is significantly larger than for classical FEM, in the practical applications classical FEM is often preferred. In \cite{farhat00} it was shown that for a given scheme that is $p$--order time--accurate on a fixed grid, satisfying the corresponding $p$--discrete space conservation law is a sufficient condition for this scheme to be at least first--order time accurate on a moving grid.

Although it seems that more work has been done in the context of the finite volumes and finite difference methods, in the context of finite elements, noticeable work has been presented in \cite{formaggia-nobile99,formaggia-nobile04,boffi-gastaldi,etienne}. \cite{etienne}. We also mention that the space conservation law is often referred to as the {\it Geometric Conservation Law} -- a phrase coined by Thomas and Lombard. It seems that various authors stick to one of these two names -- here we decided to stick to the name {\it Space Conservation Law} since it emphasizes more clearly the physics of the problem.

\subsection{Introduction of SCL}
For the numerical approximation of the equation solution, domain is {\it triangularized} into a discrete grid and the numerical solution is an array of number values attached to the grid points. Apart from numerically solving the equations of interest, two additional equations come into play when dealing with moving domains. These two additional equations pose a balance between the relevant geometric parameters -- the {\it  surface conservation law} (SCL$_s$) and the {\it volume conservation law} (SCL$_v$) -- and, as noted and discussed in \cite{zhang}, a numerical scheme which does not satisfy them shall produce additional numerical errors in the flow field. The violation of the SCL$_s$ leads to misrepresentation of the convective fluxes while violation of SCL$_v$ introduces artificial sources or sinks in otherwise conserved medium. Together, volume and surface conservation laws define the space conservation law.

Let $K_t\subset\Omega_t$ be an arbitrary control volume -- on a discrete level we are looking at cell (triangle or tetrahedra) of the triangulation ${\cal T}_h(t)$ of $\Omega(t)$ -- and denote domain velocity   at time $t$ by $\bfw(t)$. Then the time variation of the control/cell volume in terms of its boundary properties (orientations, velocities and areas) is given by
\begin{equation}\label{intgcl}
\frac{\textrm{d}}{\dt}\int_{K_t}\dbfx = \int_{\partial K_t}\bfw\cdot\bfn\dS = \int_{K_t}\div\bfw\dbfx,
\end{equation}
which is the integral statement of SCL. On the continuous level, relation (\ref{intgcl}) is trivially satisfied as long as the ALE map is regular enough. On the discrete level, however, this doesn't have to be the case and consequently artificial sources/sinks appear and (possibly) significantly influence the solution.

\subsection{SCL in finite element method}
Recall the Euler expansion formula (\ref{eulerexpansion})
\[ \frac{\partial}{\partial t}\kjac_t = \kjac_t\widehat{\div\bfw}, \]
which can be interpreted as an evolution law for the Jacobian if the domain velocity is known. It is worthy to remember that Jacobian holds the information on volume changes during the coordinate transformation.

Employing standard procedures in finite element method, the weak (conservative) form of above identity can be obtained,
\begin{equation}\label{SCLref}
\frac{\textrm{d}}{\dt}\int_{\kK}\kpsi\kjac_t\dkbfx = \int_{\kK}\kpsi\kjac_t\widehat{\div\bfw}\dkbfx,
\end{equation}
for $\kK\subset\kOmega$ and $\kpsi\in V(\kOmega)$ and, using $K_t = \kALE_t(\kK)$ to map everything onto the current configuration,
\begin{equation}\label{FEMGCL}
\frac{\textrm{d}}{\dt}\int_{K_t}\psi\dbfx = \int_{K_t}\psi\div\bfw\dbfx,
\end{equation}
where $\psi = \kpsi\circ\kALEinvt\in V(\Omega_t)$.

As mentioned in the introduction, for the technique derived in this paper the environment on referential configuration is more appropriate. Therefore, we stick to the SCL pulled back on the referent configuration i.e. with the identity (\ref{SCLref}) obtained directly from the Euler expansion formula. Applying the theory from Section~\ref{sec:cootransform}, we can obtain equivalent yet more convenient form
\begin{equation}\label{weakSCLref}
\frac{\textrm{d}}{\dt}\int_{\kK}\kpsi\kjac_t\dkbfx = \int_{\kK}\kpsi\kdiv(\kF_t\kbfw)\dkbfx.
\end{equation}
This is exactly the weak form of the differential form of SCL derived by Thomas and Lombard in \cite{thomas-lombard}. They start from the Euler expansion formula in integral form and use the metric coefficients to transform it into the form equivalent to
\begin{equation}\label{diffSCL}
\frac{\partial}{\partial t}\kjac_t - \kdiv(\kF_t\kbfw) = 0 \textrm{ in }\kOmega\times(0,T).
\end{equation}

\subsection{Conservative vs non-conservative weak formulations}
Recall for a moment the heat equation in ALE frame (\ref{strong-heat}) with the corresponding conservative weak formulation (\ref{simpleheatweak})
\[ \frac{\td}{\dt}\int_{\Omega(t)}\psi u \dbfx +\int_{\Omega(t)}\left[ \alpha\nabla\psi\cdot\nabla u - \psi\bfw\cdot\nabla u - \psi u \div\bfw -\psi f\right]\dbfx = 0. \]
Notice that an alternative, non--conservative, weak formulation is possible to obtain from (\ref{strong-heat}) by keeping the temporal derivative under the integral sign:
\begin{equation}\label{heatweaknonconservative}
\int_{\Omega(t)}\psi \frac{\partial}{\partial t}\Big|_{\kbfx} u \dbfx +\int_{\Omega(t)}\left[ \alpha\nabla\psi\cdot\nabla u - \psi\bfw\cdot\nabla u -\psi f\right]\dbfx = 0.
\end{equation}
It was argued in \cite{formaggia-nobile99} that the non--conservative formulation is able to represent the constant solution independently of the numerical time integration formula used, so it automatically satisfies the SCL. Therefore, the problematic arising from discrete SCL is not an issue in non--conservative formulation. However, the conservative discrete system maintaines  the conservation property of the original system and is therefore often preferred ahead of the non--conservative one.

\section{Formulation with vanishing discrete SCL}\label{sec:vanishSCL}
In this section the idea for ALE FEM formulation with vanishing discrete SCL is presented in details. The main objectives of the paper are given in this section.  Firstly, we show that independently on the chosen scheme for the discretization of temporal derivative, it is always possible to satisfy the SCL. Moreover, the satisfaction of SCL will be trivial since the formulation is built on the differential statement of SCL, that is the identity (\ref{diffSCL}) is in its core. Secondly, we give two different ways of calculating the grid velocity. In most approaches the calculation of the grid velocity is in tight relation with its possibility to satisfy the discrete SCL and is strongly scheme dependent, we will see that here this relation is somewhat weaker and "more decoupled" from the discretization scheme. This is thanks to the new style of the weak formulation which seems to give more freedom and possibilities for the grid velocity calculation.

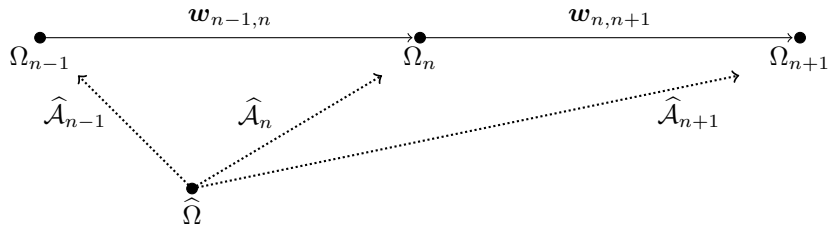
\begin{figure}[H]
\center \begin{tikzpicture}
\draw[black,fill=black] (0,0) circle (2pt) node[below] {$\Omega_{n-1}$};
\draw[black,fill=black] (5,0) circle (2pt) node[below] {$\Omega_{n}$};
\draw[black,fill=black] (10,0) circle (2pt) node[below] {$\Omega_{n+1}$};
\draw[black,fill=black] (2,-2) circle (2pt) node[below] {$\kOmega$};
\draw[black,->] (0,0)--(4.9,0);
\draw[black,->] (5,0)--(9.9,0);
\draw[black,thick,densely dotted,->] (2,-2)--(0.5,-0.5);
\draw[black,thick,densely dotted,->] (2,-2)--(4.5,-0.5);
\draw[black,thick,densely dotted,->] (2,-2)--(9.2,-0.5);

\draw[black] (1,-1) node[left] {$\kALE_{n-1}$};
\draw[black] (3.2,-1) node[left] {$\kALE_{n}$};
\draw[black] (8,-1) node[right] {$\kALE_{n+1}$};

\draw[black] (2.5,0) node[above] {$\bfw_{n-1,n}$};
\draw[black] (7.5,0) node[above] {$\bfw_{n,n+1}$};
\end{tikzpicture}
\caption{Evolution of configurations on a time interval $[t_{n-1},t_{n+1}]$.}
\label{fig-n-n+1}
\end{figure}

Before getting into details, we present the notation that will be used. Consider a set--up as that schematically shown in Figure \ref{fig-n-n+1}.
We use index $n$ to denote the functions at time $t_n$, e.g.
\begin{equation}
\begin{split}
\kALE_n & = \kALE(\cdot,t_n), \\
\kbfu_n & = \kbfu(\cdot,t_n), \\
\Omega_n & =\Omega(t_n).
\end{split}
\end{equation}
Whenever we consider a function between two time instants (functions defined piecewise), we specifically denote both of them, e.g. for the grid velocity and displacement on the time interval $[t_n,t_{n+1}]$ we write
\begin{equation}
\kbfw_{n,n+1}(\cdot,t) = \kbfw(\cdot,t_n+t)\textrm{ , } t\in[0,t_{n+1}-t_n],
\end{equation}
and
\begin{equation}
\kbfu_{n,n+1}(\cdot,t) = \kbfu(\cdot,t_n+t) \textrm{ , }t\in[0,t_{n+1}-t_n],
\end{equation}
respectively. Note that the relation between $\kbfu_{n,n+1}$ and $\kbfw_{n,n+1}$ is
\begin{equation}
\kbfu_{n,n+1}(\cdot,t) = \kbfu_n(\cdot) + \int_{0}^{t_n+t}\kbfw_{n,n+1}(\cdot,s)\ds\textrm{ , }t\in[0,t_{n+1}-t_n].
\end{equation}
As mentioned before, using the {\it "hat"} symbol or dropping it denotes which configuration a function is defined on -- {\it "hat"} for the functions on the referential configuration while dropping it for the functions on the physical configurations.

Now, in the spirit of the above introduced notation, we can write
\begin{equation}
\begin{split}
\kALE(\kbfx,t_{n}+t) = \kALE_{n,n+1}(\kbfx,t) & = \kbfx + \kbfu_{n,n+1}(\kbfx,t)\textrm{ , }t\in[0,t_{n+1}-t_n].
\end{split}
\end{equation}

In the next two subsections we give two possibilities for the grid velocity calculation -- the first being the classical one used in most approaches, while the second is a bit more advanced and physically very reasonable. We also show that for these choices the SCL identities -- weak form (\ref{weakSCLref}) and differential form (\ref{diffSCL}) -- are trivially satisfied when handled properly.

\subsection{Grid velocity calculation and satisfaction of SCL}\label{sec:grid-vel-dc}

Denote the time step between two time instants by $\Delta t_n^{n+1}=t_{n+1}-t_n$. Then the most widely used method for the grid velocity calculation states: knowing the position of the grid node indexed by $i$ at times $t_n$ and $t_{n+1}$, i.e. the values $\bfx_i^n=\kbfx_i+\kbfu_n(\kbfx_i)$ and $\bfx_i^{n+1}=\kbfx_i+\kbfu_{n+1}(\kbfx_i)$, the grid velocity can be defined as piecewise constant in time on $[t_n,t_{n+1}]$ by
\begin{equation}\label{const-grid-vel}
\kbfw_{n,n+1}(\kbfx_i,t) = \frac{\kbfu_{n+1}(\kbfx_i)-\kbfu_{n}(\kbfx_i)}{\Delta t_n^{n+1}}\textrm{ , }t\in[0,\Delta t_n^{n+1}].
\end{equation}
In that case
\begin{equation}
\kbfu_{n,n+1}(\cdot, t) = \kbfu_n + t\kbfw_{n,n+1}\textrm{ , }t\in[0,\Delta t_n^{n+1}]
\end{equation}
and 
\begin{equation}
\kALE_{n+1}(\kbfx) = \kbfx + \kbfu_n(\kbfx) + (\Delta t_n^{n+1})\kbfw_{n,n+1}(\kbfx).
\end{equation}
The evolution between two configurations $\Omega_{n}$ and $\Omega_{n+1}$ is given by
\begin{equation}
\Omega(t_n+t) = \kALE_{n,n+1}(\kOmega,t) \textrm{ , }t\in[0,\Delta t_n^{n+1}]
\end{equation}
where 
\begin{equation}
\kALE_{n,n+1}(\cdot,t) = \kbfx + \kbfu_n(\kbfx) + t\kbfw_{n,n+1}(\kbfx)\textrm{ , }t\in[0,\Delta t_n^{n+1}],
\end{equation}
i.e. $\kALE(\cdot,t+t_n) = \kALE_{n,n+1}(\cdot,t)$, $t\in[0,\Delta t_n^{n+1}]$.

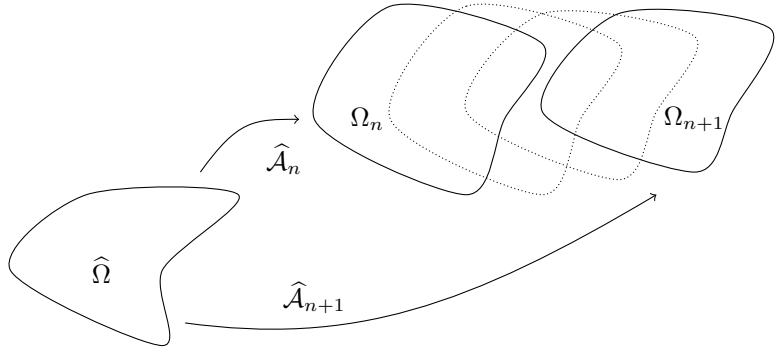
\begin{figure}[H]
\center \begin{tikzpicture}
\draw [black] plot [smooth cycle] coordinates {(0,0) (1,1) (3,1) (2,0) (2,-1)};

\draw [black] plot [smooth cycle] coordinates {(4,2) (5,3.5) (7,3) (6.5,2) (6,1)};
\draw [black] plot [smooth cycle] coordinates {(7,2.1) (8,3.4) (10,3.2) (9.5,2.1) (9,1.3)};

\draw [black,densely dotted] plot [smooth cycle] coordinates {(5,2) (6,3.5) (8,3) (7.5,2) (7,1)};
\draw [black,densely dotted] plot [smooth cycle] coordinates {(6,2.1) (7,3.4) (9,3.1) (8.5,2.1) (8,1.2)};

\draw[black] (1.2,0) node[] {$\kOmega$};
\draw[black] (4.7,2) node[] {$\Omega_n$};
\draw[black] (9,2) node[] {$\Omega_{n+1}$};

\draw[black,->] (2.5,1.3).. controls (3,2) and (3.2,2) ..(3.8,2);
\draw[black] (3.6,1.8) node[below] {$\kALE_n$};

\draw[black,->] (2.3,-0.7).. controls (4.5,-1) and (6,-0.5) ..(8.5,1);
\draw[black] (4,-0.7) node[above] {$\kALE_{n+1}$};

\end{tikzpicture}
\caption{Sketch of a configuration evolution between $t_n$ and $t_{n+1}$.}
\end{figure}

Now consider the differential form of SCL identity (\ref{diffSCL}),
\[ \frac{\partial}{\partial t}\kjac_t - \kdiv(\kF_t\kbfw) = 0 \textrm{ in }\kOmega\times(0,T), \]
and let us take a more detailed look on what is happening during the temporal discretization. For the sake of notational simplicity, assume that $t_{n+1}-t_n=\Delta t$, $\forall n$.
Discretized counterpart of $\frac{\partial}{\partial t}\kjac_t$ takes the form of
\[ \frac{\partial}{\partial t}\kjac_t \approx \frac{\kjac_{n+1} - \kjac_n}{\Delta t} \]
which is obtained from some integration quadrature formula. In case of implicit Euler formula, for example, using the fact that $\kF_t$ and $\kbfw$ are defined piecewise, we have
\begin{equation}\label{IEintegrationformula}
\begin{split}
\kjac_{n+1} - \kjac_n = \int_{t_n}^{t_{n+1}} \frac{\partial}{\partial t}\kjac_t \dt  & = \int_{t_n}^{t_{n+1}} \kdiv(\kF_{n,n+1}(t)\kbfw_{n,n+1}(t)) \dt\\
& \approx \Delta t \kdiv(\kF_{n,n+1}(t_{n+1})\kbfw_{n,n+1}(t_{n+1})).
\end{split}
\end{equation}
Taking into account that $ \kF_{n,n+1}(t)$ and $\kbfw_{n,n+1}(\cdot,t)$ are polynomials in $t$ variable, so is the $\kF_{n,n+1}(t)\kbfw_{n,n+1}(t)$. Then we can see that the step in approximating the time integral in (\ref{IEintegrationformula}) is actually unnecessary since the integral can formally interchange with the divergence $\kdiv$ operator and the integral of polynomial can be evaluated exactly. Therefore, we can write 
\begin{equation}
\kjac_{n+1} - \kjac_n =  \kdiv\left[\int_{t_n}^{t_{n+1}}\kF_{n,n+1}(t)\kbfw_{n,n+1}(t)\dt\right]
\end{equation}
and the discrete SCL vanishes. The same procedure can be applied to the weak form of SCL, the identity (\ref{weakSCLref}):
\[ \frac{\textrm{d}}{\dt}\int_{\kK}\kpsi\kjac_t\dkbfx = \int_{\kK}\kpsi\kdiv(\kF_t\kbfw)\dkbfx. \] Integrating this identity from $t_n$ to $t_{n+1}$ and noticing that we can formally change the order of integration (since referential control volumes $\kK$ are independent of time) and that the test function $\kpsi$ is time independent (so behaving like a constant with respect to $\int_{t_n}^{t_{n+1}}\dt$) we obtain
\begin{equation}
\int_{\kK}\kpsi(\kjac_{n+1}-\kjac_{n})\dkbfx = \int_{\kK}\kpsi\kdiv\left[\int_{t_n}^{t_{n+1}}\kF_{n,n+1}(t)\kbfw_{n,n+1}(t)\dt\right]\dkbfx.
\end{equation}
Again, since polynomials can be integrated exactly, the discrete SCL vanishes and is trivially satisfied in this type of formulation.

\subsection{Grid velocity continuous in time}

Consider for a moment the following situation: grid node $\bfx$ moves with velocity $\bfw_{n-1,n}\neq 0$ from its position at time $t_{n-1}$ to a new position at time $t_{n}$ during the time step of length $\Delta t=t_n-t_{n-1}$. Then the new position at time $t_{n+1}$ can be somehow (depending on the problem of interest) obtained. Assume for a moment that the position at time $t_{n+1}$ is equal to the position at time $t_n$, that is $\bfx^{n+1}=\bfx^{n}$. Using the method introduced in Section~\ref{sec:grid-vel-dc}, we would find that \[\bfw_{n,n+1}(\bfx,t)=0.\] However, a physically more reasonable explanation would be that the time step $\Delta t$ is too large for the numerics to "catch" the entire trajectory of node $\bfx$ (if node stops instantly then its "trajectory momentum" is violated). Since the node traveled to position $\bfx^n$ with velocity $\bfw_{n-1,n}(\cdot,t)$, $\bfw_{n-1,n}(\cdot,t_n)\neq 0$ , it is reasonable to assume that its velocity on the time interval $[t_n,t_{n+1}]$, $\bfw_{n,n+1}(\cdot,t)$, is continuous at $t_n$ i.e. $\neq 0$. Thus, the trajectory of the node $\bfx$ is "a closed loop" starting and ending at $\bfx^n=\bfx^{n+1}$. Motivated by this example, we propose an alternative approach to calculate the grid velocity continuous in time once the positions of the nodes are known.

Let $\bfw_{n-1,n}(\cdot,t)$ on $[t_{n-1},t_n]$ be known, and assume $\Omega_{n+1}$ has been found. Then we define the grid velocity on interval $[t_n,t_{n+1}]$ by
\begin{equation}\label{con-grid-vel}
\bfw(\cdot,t) = (t-t_n)\bfomega_{n,n+1}+ \bfw(\cdot,t_n)\textrm{ , }t\in[t_n,t_{n+1}].
\end{equation} Its "piecewise counterpart" is then expressed as
\begin{equation}
\bfw_{n,n+1}(\cdot,t)=t\bfomega_{n,n+1} + \bfw(\cdot,t_n)\textrm{ , }t\in[0,t_{n+1}-t_n],
\end{equation}
where $\bfomega_{n,n+1}\in\mR^2$  is a constant (in time) appropriately chosen vector field (defined bellow). The first consequence to notice is that the grid velocity is now continuous, i.e.
\begin{equation}
\bfw_{n,n+1}(\cdot,0) = \bfw_{n-1,n}(\cdot,\Delta t_{n-1}^n).
\end{equation}
Next we have to determine the constant field $\bfomega_{n,n+1}$. Clearly, we have to demand 
\begin{equation}
\int_0^{t_{n+1}-t_n}\bfw_{n,n+1}(\cdot,t)\dt = \kbfu_{n+1}(\cdot)-\kbfu_n(\cdot).
\end{equation} 
By a straightforward calculation, it follows
\begin{equation}
\bfomega_{n,n+1}(\cdot) = \frac{2}{(t_{n+1}-t_n)^2}\big[ \kbfu_{n+1}(\cdot)-\kbfu_n(\cdot) - (t_{n+1}-t_n)\bfw(\cdot,t_n)\big]. 
\end{equation}
Then using the relation
\[ \kbfu_{n,n+1}(\cdot,t) = \kbfu_n(\cdot)+\int_0^{t}\kbfw_{n,n+1}(\cdot,s)\ds\textrm{ , }t\in[0,t_{n+1}-t_n] \]
we interpolate the displacement on whole $[t_n,t_{n+1}]$
\begin{equation}
\kbfu(\cdot,t_n+t)=\kbfu_n(\cdot)+ t\kbfw_{n-1,n}(\cdot,t_n) + \frac{t^2}{2}\kbfomega_{n,n+1}(\cdot) \textrm{ , }t\in[0,t_{n+1}-t_n].
\end{equation}
Now analogously as in Section~\ref{sec:grid-vel-dc} we deduce that $\kF_{n,n+1}(t)\kbfw_{n,n+1}(t)$ is piecewise polynomial in time so the argument for vanishing discrete SCL stays the same.

\section{Discretization schemes}\label{sec:schemes}
In this section we aim to show how to properly handle some classic discretization schemes in order for the SCL to be preserved. We deal with the {\it implicit (backward) Euler} scheme, the {\it Crank--Nicolson} scheme, and the {\it backward differentiation formulas}, {\it BDF2} and {\it BDF3}. For simplicity, all of these schemes shall be illustrated on the heat equation. Generalizations to {\it convection--diffusion} or {\it Navier--Stokes} equations is straightforward since handling the terms with grid velocity stays the same and only these terms play a role in problematics in the context of SCL. 

Recall the heat equation in the weak form pulled to the referent configuration:
\begin{equation}\label{main-eq}
\begin{split}
0=\frac{\td}{\dt} & \int_{\kOmega}  \kpsi\ku\kjac_t  \dkbfx \\
+ & \int_{\kOmega}\left[ \alpha\frac{1}{\kjac_t}\kF_t\kF_t^T\knabla\kpsi\cdot\knabla\ku - \kpsi\kF_t\kbfw\cdot\knabla\ku  - \kpsi\ku\kdiv(\kF_t\kbfw) - \kpsi\kf\kjac_t\right]\dkbfx,
\end{split}
\end{equation}
with the initial and boundary conditions
\begin{equation}
\begin{split}
\ku(0) & = \ku_0 \textrm{ in }\kOmega, \\
\ku & = \ku_D\textrm{ on }\partial\kOmega\times(0,T).
\end{split}
\end{equation}

Before getting into the details with discretization schemes, consider for a moment the diffusion term
\[ \alpha\frac{1}{\kjac_t}\kF_t\kF_t^T\knabla\kpsi\cdot\knabla\ku \] in the above weak formulation. Due to $\frac{1}{\kjac_t}$ we see that rational function in $t$ will appear, which generally cannot be integrated exactly. Still, we can employ some quadrature formula for integration (e.g. {\it Simpsons rule}), but since it does not play a significant role in SCL problematics, we simply approximate it implicitly for the moment, i.e.
\begin{equation}
\begin{split}
\int_{t_n}^{t_{n+1}}\int_{\kOmega} & \alpha\frac{1}{\kjac_t}\kF_t\kF_t^T\knabla\kpsi\cdot\knabla\ku\dkbfx\dt \\
& \approx \Delta t \int_{\kOmega} \alpha\frac{1}{\kjac_{n,n+1}(\Delta t)}\kF_{n,n+1}(\Delta t)\kF_{n,n+1}^T(\Delta t)\knabla\kpsi\cdot\knabla\ku_{n+1}\dkbfx 
\end{split}
\end{equation}
with $\Delta t = t_{n+1}-t_n$. The same we apply to the source term involving $\kf$. For the compact notation we denote
\begin{equation}
d_{n,n+1}(\ku_k,\kpsi) = \Delta t \int_{\kOmega} \alpha\frac{1}{\kjac_{n,n+1}(\Delta t)}\kF_{n,n+1}(\Delta t)\kF_{n,n+1}^T(\Delta t)\knabla\kpsi\cdot\knabla\ku_k\dkbfx
\end{equation}
where $\ku_k$ is the approximation of $\ku(t_k)$ after the temporal discretization, and $k$ is typically $n$ or $n+1$, depending on the chosen discretization schemes. Similarly, we denote
\begin{equation}
b_{n,n+1}(\kf_k,\kpsi) = \Delta t\int_{\kOmega}\kpsi\kf_k\kjac_{n,n+1}(\Delta t)\dkbfx
\end{equation}
with $\kf_k=\kf(t_k)$. For terms involving grid velocity and operators related to ALE map, we denote
\begin{equation}
\begin{split}
{\cal M}_{n,n+1}(\ku_k,\kpsi) = & \int_{\kOmega} \kpsi  \left[\int_{0}^{\Delta t_n^{n+1}}\kF_{n,n+1}(t)\kbfw_{n,n+1}(t)\dt\right]\cdot\knabla\ku_k\dkbfx \\
& + \int_{\kOmega}\kpsi\ku_k\kdiv\left[\int_{0}^{\Delta t_n^{n+1}}\kF_{n,n+1}(t)\kbfw_{n,n+1}(t)\dt\right]\dkbfx
\end{split}
\end{equation}
The general procedure of the new approach applied to the classical schemes is to modify the terms involving the grid velocity through ${\cal M}_{n,n+1}$ operator, while the rest remains in the spirit of the original scheme.

\subsection{Implicit Euler method}
The implicit Euler scheme is obtained by integrating (\ref{main-eq}) from $t_n$ to $t_{n+1}$ and approximating temporal integrals implicitly, i.e. unknown under the integral sign is taken at time $t_{n+1}$. Thus, by implicit Euler scheme we obtain (in semidiscrete version)
\begin{equation}
\begin{split}
\int_{\kOmega}\kpsi\ku_{n+1} & \kjac_{n+1}\dkbfx -  \int_{\kOmega}\kpsi\ku_{n}\kjac_{n}\dkbfx \\
& +d_{n,n+1}(\ku_{n+1},\kpsi) - b_{n,n+1}(\kf_{n+1},\kpsi) - {\cal M}_{n,n+1}(\ku_{n+1},\kpsi)=0
\end{split}
\end{equation}
\begin{figure}[H]
\center \begin{tikzpicture}
\draw[black,fill=black] (0,0) circle (2pt) node[below] {$\Omega_{n-1}$};
\draw[black,fill=black] (5,0) circle (2pt) node[below] {$\Omega_{n}$};
\draw[black,fill=black] (10,0) circle (2pt) node[below] {$\Omega_{n+1}$};
\draw[black,densely dotted,->] (0,0)--(4.9,0);
\draw[black,->] (5,0)--(9.9,0);

\draw[black,densely dotted] (0,1)--(4.9,1);
\draw[black,densely dotted] (5,2)--(9.9,2);

\draw[black] (2.5,0) node[below] {$\bfw_{n-1,n}(t)$};
\draw[black] (2.5,-0.5) node[below] {$\jac_{n-1,n}(t)$};

\draw[black] (0,-0.5) node[below] {$\jac_{n-1}$};
\draw[black] (5,-0.5) node[below] {$\jac_{n}$};
\draw[black] (10,-0.5) node[below] {$\jac_{n+1}$};

\draw[black] (7.5,0) node[below] {$\bfw_{n,n+1}(t)$};
\draw[black,fill=black] (5,1) circle (2pt);
\draw[black] (0,1) circle (2pt);
\draw[black] (7.5,-0.5) node[below] {$\jac_{n,n+1}(t)$};
\draw[black,fill=black] (10,2) circle (2pt);
\draw[black] (5,2) circle (2pt);

\draw[black] (2.5,1) node[above] {$u_n$};
\draw[black] (7.5,2) node[above] {$u_{n+1}$};
\end{tikzpicture}
\caption{Sketch of the implicit Euler method on $[t_{n-1},t_{n+1}]$. When calculating $u_{n+1}$ we see that all the action is happening on $[t_n,t_{n+1}]$ and test/basis functions involved in SCL "carry the same weight" so there is no violation of SCL.}
\label{IEsketch}
\end{figure}
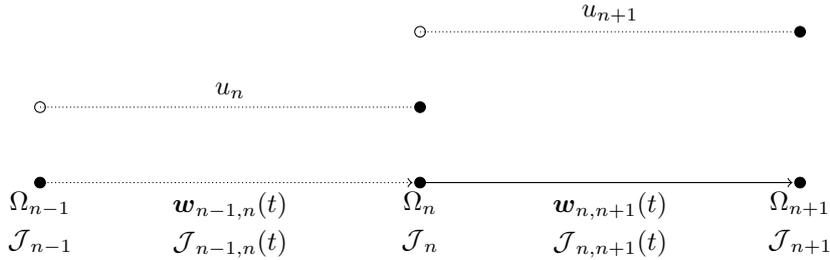

\subsection{Crank--Nicolson method}
Crank--Nicolson scheme is obtained by integrating (\ref{main-eq}) from $t_n$ to $t_{n+1}$ and approximating temporal integrals by trapezoidal rule. Thus, by Crank--Nicolson scheme we obtain (in semidiscrete version)
\begin{equation}
\begin{split}
\int_{\kOmega}\kpsi & \ku_{n+1}\kjac_{n+1}\dkbfx -  \int_{\kOmega}\kpsi  \ku_{n}\kjac_{n}\dkbfx \\
& +\frac{1}{2} \left[ d_{n,n+1}(\ku_{n+1},\kpsi)+d_{n,n+1}(\ku_{n},\kpsi)\right]  - \frac{1}{2}\left[ b_{n,n+1}(\kf_{n+1},\kpsi) + b_{n,n+1}(\kf_{n},\kpsi) \right] \\
 &- \frac{1}{2}\left[ {\cal M}_{n,n+1}(\ku_{n+1},\kpsi) + {\cal M}_{n,n+1}(\ku_{n},\kpsi)\right] = 0.
\end{split}
\end{equation}
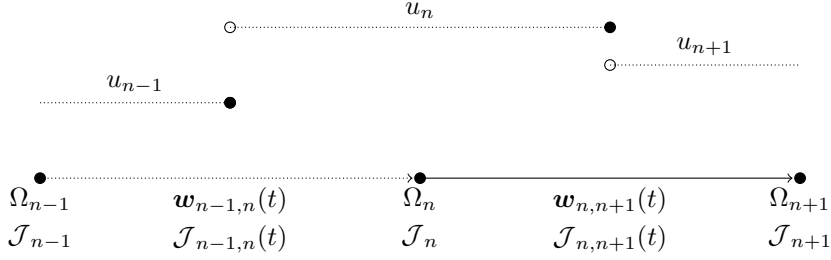
\begin{figure}[H]
\center \begin{tikzpicture}
\draw[black,fill=black] (0,0) circle (2pt) node[below] {$\Omega_{n-1}$};
\draw[black,fill=black] (5,0) circle (2pt) node[below] {$\Omega_{n}$};
\draw[black,fill=black] (10,0) circle (2pt) node[below] {$\Omega_{n+1}$};
\draw[black,densely dotted,->] (0,0)--(4.9,0);
\draw[black,->] (5,0)--(9.9,0);

\draw[black,densely dotted] (0,1)--(2.5,1);
\draw[black,densely dotted] (2.5,2)--(7.5,2);
\draw[black,densely dotted] (7.5,1.5)--(10,1.5);

\draw[black] (2.5,0) node[below] {$\bfw_{n-1,n}(t)$};
\draw[black] (2.5,-0.5) node[below] {$\jac_{n-1,n}(t)$};

\draw[black] (0,-0.5) node[below] {$\jac_{n-1}$};
\draw[black] (5,-0.5) node[below] {$\jac_{n}$};
\draw[black] (10,-0.5) node[below] {$\jac_{n+1}$};

\draw[black] (7.5,0) node[below] {$\bfw_{n,n+1}(t)$};
\draw[black,fill=black] (2.5,1) circle (2pt);
\draw[black] (2.5,2) circle (2pt);
\draw[black] (2.5,1) circle (2pt);
\draw[black] (7.5,-0.5) node[below] {$\jac_{n,n+1}(t)$};
\draw[black] (7.5,1.5) circle (2pt);
\draw[black,fill=black] (7.5,2) circle (2pt);

\draw[black] (1.25,1) node[above] {$u_{n-1}$};
\draw[black] (5,2) node[above] {$u_{n}$};
\draw[black] (8.75,1.5) node[above] {$u_{n+1}$};
\end{tikzpicture}
\caption{Sketch of the Crank--Nicolson method on $[t_{n-1},t_{n+1}]$. When calculating $u_{n+1}$ we see that all the action is happening on $[t_n,t_{n+1}]$ and test/basis functions involved in SCL "carry the same weight" so there is no violating SCL.}
\label{CNsketch}
\end{figure}
\subsection{Backward differentiation formula -- BDF}
With the backward differentiation formulas, the situation is a little bit different, primarily because, in contrast to the  previous two methods, it is based on a differentiation instead of integration. Consider first the BDF2 method and the ordinary differential equation (ODE) in the form of $y'(t)=f(t,y(t))$. Then the ODE discretized by the BDF2 takes form
\begin{equation}\label{genBDF2}
\frac{3}{2}y^{n+1} - 2 y^n + \frac{1}{2}y^{n-1} = \Delta t f(t_{n+1},y^{n+1}).
\end{equation}
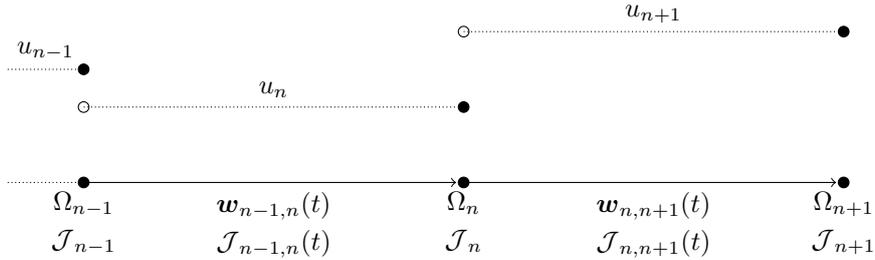
\begin{figure}[H]
\center \begin{tikzpicture}
\draw[black,fill=black] (0,0) circle (2pt) node[below] {$\Omega_{n-1}$};
\draw[black,fill=black] (5,0) circle (2pt) node[below] {$\Omega_{n}$};
\draw[black,fill=black] (10,0) circle (2pt) node[below] {$\Omega_{n+1}$};
\draw[black,->] (0,0)--(4.9,0);
\draw[black,->] (5,0)--(9.9,0);

\draw[black,densely dotted] (0,1)--(4.9,1);
\draw[black,densely dotted] (5,2)--(9.9,2);

\draw[black] (2.5,0) node[below] {$\bfw_{n-1,n}(t)$};
\draw[black] (2.5,-0.5) node[below] {$\jac_{n-1,n}(t)$};

\draw[black] (0,-0.5) node[below] {$\jac_{n-1}$};
\draw[black] (5,-0.5) node[below] {$\jac_{n}$};
\draw[black] (10,-0.5) node[below] {$\jac_{n+1}$};

\draw[black] (7.5,0) node[below] {$\bfw_{n,n+1}(t)$};
\draw[black,fill=black] (5,1) circle (2pt);
\draw[black] (0,1) circle (2pt);
\draw[black] (7.5,-0.5) node[below] {$\jac_{n,n+1}(t)$};
\draw[black,fill=black] (10,2) circle (2pt);
\draw[black] (5,2) circle (2pt);

\draw[black,fill=black] (0,1.5) circle (2pt);
\draw[black,densely dotted] (-1,1.5)--(0,1.5);
\draw[black,densely dotted] (-1,0)--(0,0);
\draw[black] (-0.5,1.5) node[above] {$u_{n-1}$};

\draw[black] (2.5,1) node[above] {$u_n$};
\draw[black] (7.5,2) node[above] {$u_{n+1}$};
\end{tikzpicture}
\caption{Sketch of the BDF2 method on $[t_{n-1},t_{n+1}]$. When calculating $u_{n+1}$ we see that all the action is happening on $[t_{n-1},t_{n+1}]$ so, in order not to violate the SCL, the whole evolution of configuration on $[t_{n-1},t_{n+1}]$ has to be taken into account.}
\label{BDF2sketch}
\end{figure}
Rearranging the left hand side in general form of BDF2 we can obtain
\begin{equation}
\frac{3}{2}y^{n+1} - 2 y^n + \frac{1}{2}y^{n-1} = \frac{3}{2}(y^{n+1}-y^n) - \frac{1}{2}(y^n-y^{n-1}).
\end{equation}
In the context of the time--dependent functions defined on the time--dependent domains, when we talk about a function at time $t$ we actually have in mind the pair of function and its domain at time $t$. Applying the BDF2 scheme for the time discretization of our model problem (\ref{main-eq}) we obtain
\begin{equation}
\begin{split}
\frac{3}{2}\left(\int_{\kOmega}\kpsi\ku_{n+1}\kjac_{n+1}\dkbfx  - \int_{\kOmega}\kpsi\right.&\left.\ku_{n}\kjac_{n}\dkbfx \right) \\ - & \frac{1}{2}\left(\int_{\kOmega}\kpsi\ku_{n}\kjac_{n}\dkbfx - \int_{\kOmega}\kpsi\ku_{n-1}\kjac_{n-1}\dkbfx \right)
\end{split}
\end{equation} 
from which it is clear that the action is happening on $[t_{n-1},t_n]$ and $[t_{n},t_{n+1}]$ with "weights" $\frac{3}{2}$ and $-\frac{1}{2}$, respectively. Therefore, in order for the SCL not to be violated, the full evolution with the respective weights has to be considered in (\ref{genBDF2}) on the right hand side. So the (modified) BDF2 scheme which does not violate the SCL states
\begin{equation}
\begin{split}
\frac{3}{2}\int_{\kOmega}\kpsi\ku_{n+1}\kjac_{n+1}\dkbfx & - 2\int_{\kOmega}\kpsi\ku_{n}\kjac_{n}\dkbfx + \frac{1}{2}\int_{\kOmega}\kpsi\ku_{n-1}\kjac_{n-1}\dkbfx\\
& +d_{n,n+1}(\ku_{n+1},\kpsi) - b_{n,n+1}(\kf_{n+1},\kpsi) \\
& -\frac{3}{2} {\cal M}_{n,n+1}(\ku_{n+1},\kpsi) + \frac{1}{2}{\cal M}_{n-1,n}(\ku_{n+1},\kpsi)=0.
\end{split}
\end{equation}

The procedure for BDF3 is analogous: starting from the general form
\begin{equation}\label{genBDF3}
\frac{11}{6}y^{n+1} - 3 y^n + \frac{3}{2}y^{n-1} - \frac{1}{3}y^{n-2} = \Delta t f(t_{n+1},y^{n+1})
\end{equation}
we can rewrite it as
\begin{equation}
\frac{11}{6}y^{n+1} - 3 y^n + \frac{3}{2}y^{n-1} - \frac{1}{3}y^{n-2} = \frac{11}{6}(y^{n+1}-y^n) - \frac{7}{6}(y^n-y^{n-1}) + \frac{1}{3}(y^{n-1}-y^{n-2})
\end{equation}
in order to get the final (modified) BDF3 formulation
\begin{equation}
\begin{split}
\int_{\kOmega}\left( \right. & \left. \frac{11}{6}  \kpsi  \ku_{n+1}  \kjac_{n+1}  - 3  \kpsi\ku_{n}\kjac_{n}  + \frac{3}{2}\kpsi\ku_{n-1}\kjac_{n-1} - \frac{1}{3}\kpsi\ku_{n-2}\kjac_{n-2}\right)\dkbfx\\
& +d_{n,n+1}(\ku_{n+1},\kpsi) - b_{n,n+1}(\kf_{n+1},\kpsi) \\
& -\frac{11}{6} {\cal M}_{n,n+1}(\ku_{n+1},\kpsi) + \frac{7}{6}{\cal M}_{n-1,n}(\ku_{n+1},\kpsi) - \frac{1}{3}{\cal M}_{n-2,n-1}(\ku_{n+1},\kpsi)\\
& =0.
\end{split}
\end{equation}

\section{Spatial discretization}\label{sec:fem}
For spatial discretization, finite element method has been employed. For details we refer to \cite{ciarlet,quarteroni} and  the references therein.  We mention that a good summary on construction of the finite element spaces on time--dependent domains is given in \cite{formaggia-nobile99}. Here we only summarize a few honorable mentions on the selection of the finite element spaces for ALE map and the related fields naturally arising from it (e.g. grid displacement, velocity, Jacobian of ALE map). It is important to notice that once the finite element space for ALE map $\kALE_t$ is chosen, its selection (partially) dictates the selections of finite element spaces for the fields derived from the ALE map.

Assume that the grid is triangularized into a finite number of triangles (2d) or tetrahedra (3d) with the straight edges/faces. In order to preserve the straight edges/faces during the grid motion, the displacement has to be chosen as piecewise (per triangle) linear polynomial in space, that is
\begin{equation}
\kALE_{h,t},\kbfu_{h,t},\kbfw_{h,t}\in[\mathbb{P}_1]^d.
\end{equation}
The index $h$ denotes the discretized (in space) counterpart of the function. In that case \[\kALEF_{h,t}\in[\mathbb{P}_0]^{d\times d}\] so $\kjac_{h,t}\in\mathbb{P}_0$ and $\kF_{h,t}\in[\mathbb{P}_0]^{d\times d}$. Then, since piecewise constant function times piecewise first order polynomial continuous function is generally piecewise first order polynomial but is discontinuous, we have \[\kF_{h,t}\kbfw_{h,t}\in[\mathbb{P}_{1,dc}]^d,\] $\mathbb{P}_{1,dc}$ denoting the space of piecewise first order polynomials that are discontinuous over the edges/faces of elements.

From the above discussion it is clear that we do not have  complete freedom in choosing finite element spaces for ALE map. If, say, we were to choose $\mathbb{P}_1$ space for all the fields mentioned above, most softwares do an immediate interpolation thus producing (an additional) error in the spatial discretization. So if one is not careful, the violation of SCL may occur even though we are dealing with the time integrals exactly.

\section{Numerical validation}\label{sec:numerics}
In this section we perform a numerical validation of the proposed approach. For the first and second order schemes we compare our results with the benchmark problems set--up in \cite{formaggia-nobile04,boffi-gastaldi}. First the stability for different schemes is tested, and then the accuracy is assessed. For the benchmark problems in original papers \cite{formaggia-nobile04,boffi-gastaldi} piecewise linear finite elements are employed and we stick to the same selection. 
\subsection{Stability}\label{sec:stability}
The following problem is considered:
\begin{equation}
\begin{split}
\partial_t u - 0.01\Delta u & = 0 \textrm{ in }\Omega(t)\times(0,T)\\
u & = 0\textrm{ on }\partial\Omega(t)\times(0,T)\\
u(0) & = 1600x(1-x)y(1-y)\textrm{ in }\Omega(0)
\end{split}
\end{equation}
with $\Omega(0)=[0,1]^2$. The prescribed ALE map is given below
\begin{equation}
\kALE_t(\kbfx) = (2-\cos 20\pi t)\kbfx\textrm{ in }\kOmega=\Omega(0).
\end{equation}
The time interval of interest is $[0,0.4]$ which corresponds to four periods of oscillations of the domain.

In \cite{formaggia-nobile04,boffi-gastaldi} they use the {\it Gronwall lemma} to show that the norm $\Vert u(t)\Vert_{L^2(\Omega(t))}$ decreases with $t$. Therefore, for a stable discretization, the same decreasing trend should be expected for the discrete solution. 

\begin{figure}
\center \includegraphics[scale=0.14]{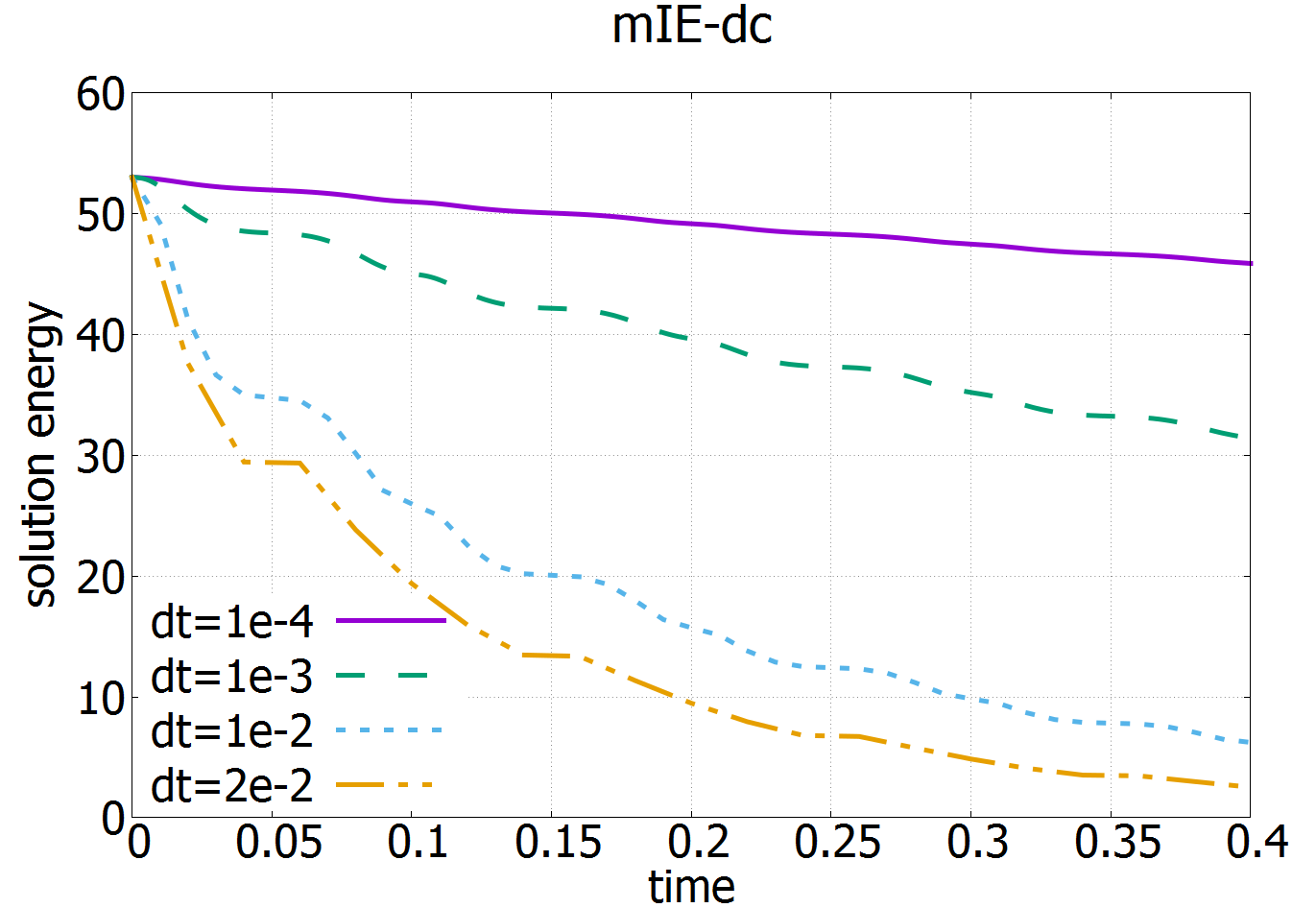}\includegraphics[scale=0.14]{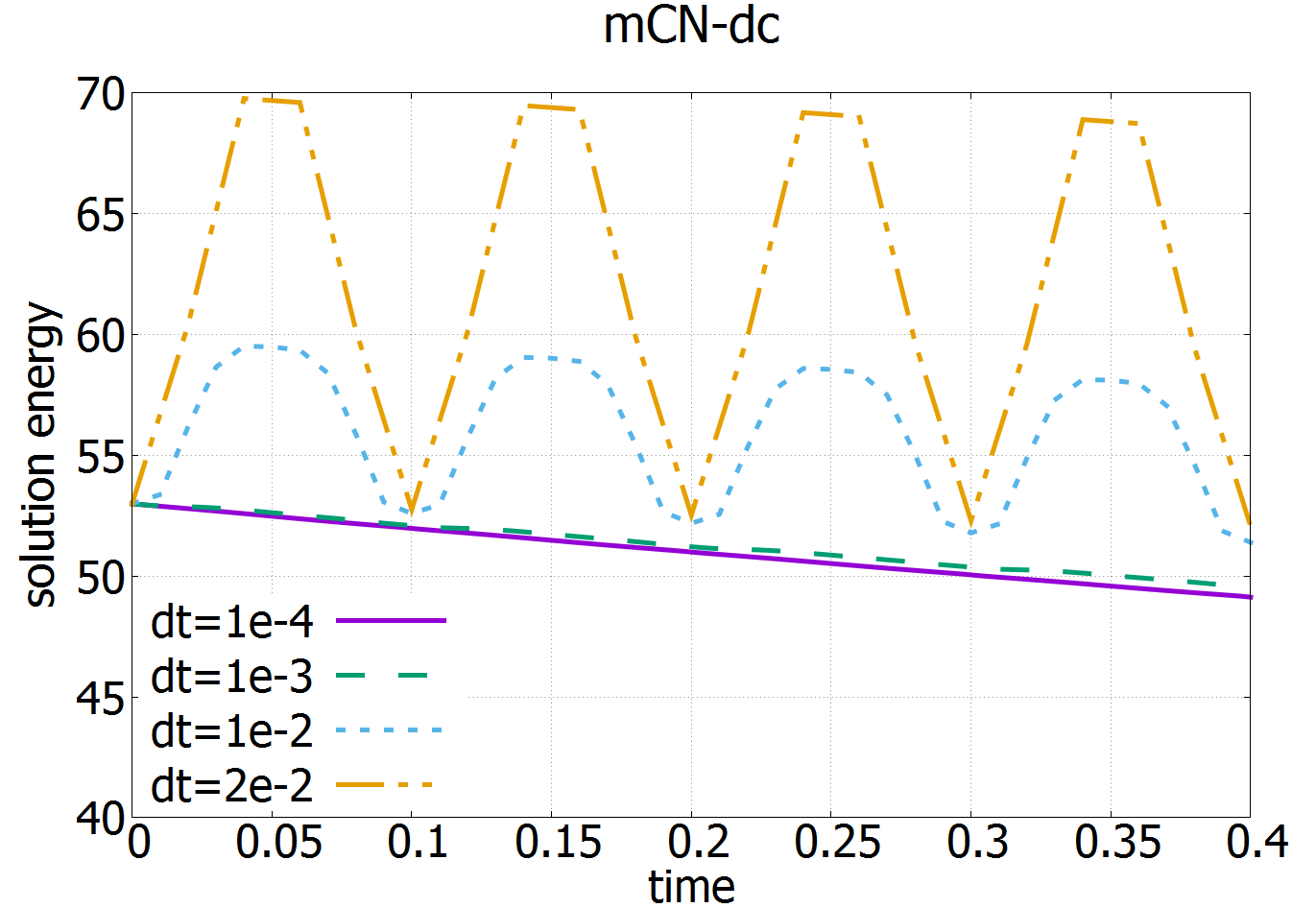}

\vspace{0.5cm}

\center \includegraphics[scale=0.14]{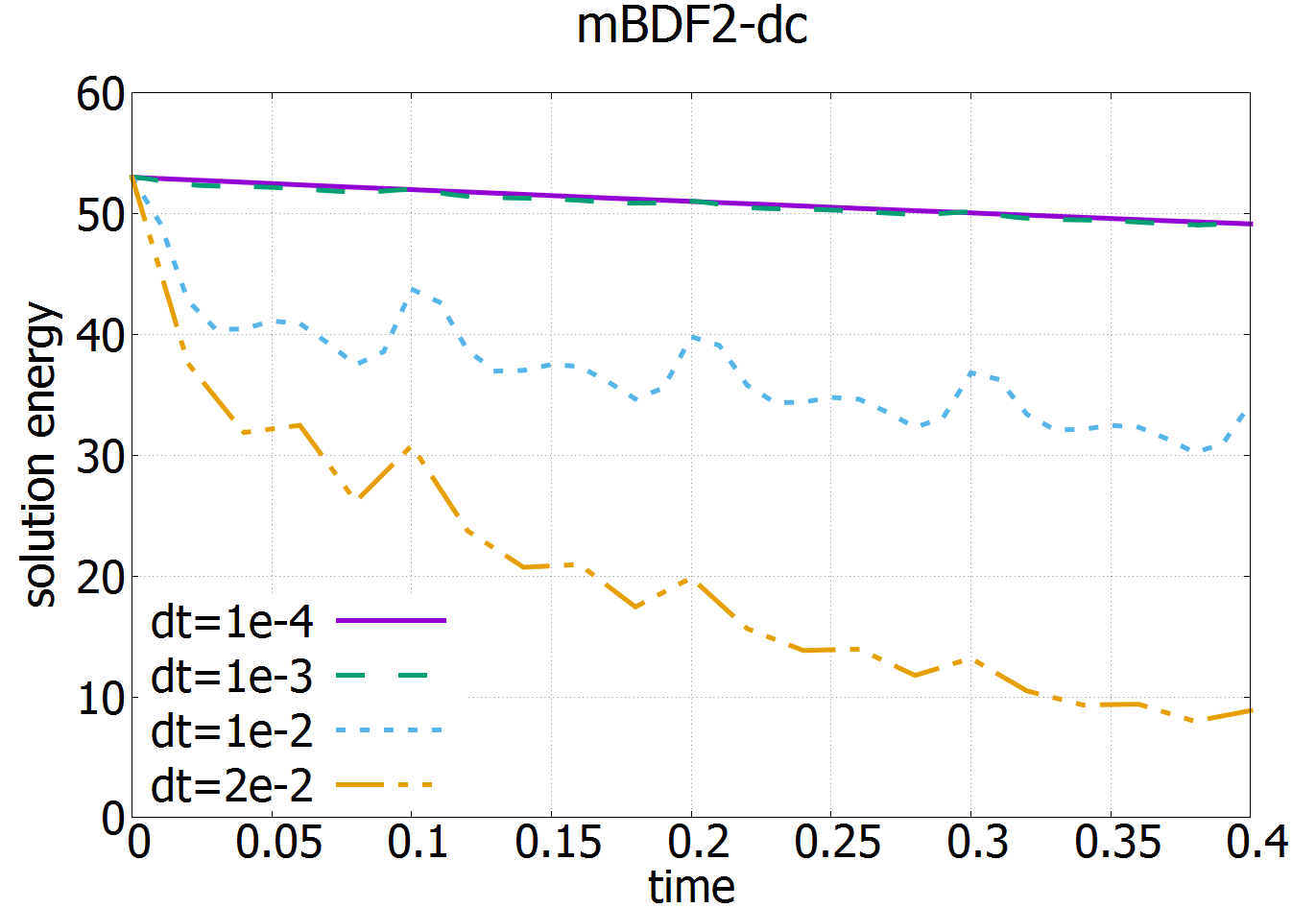}\includegraphics[scale=0.14]{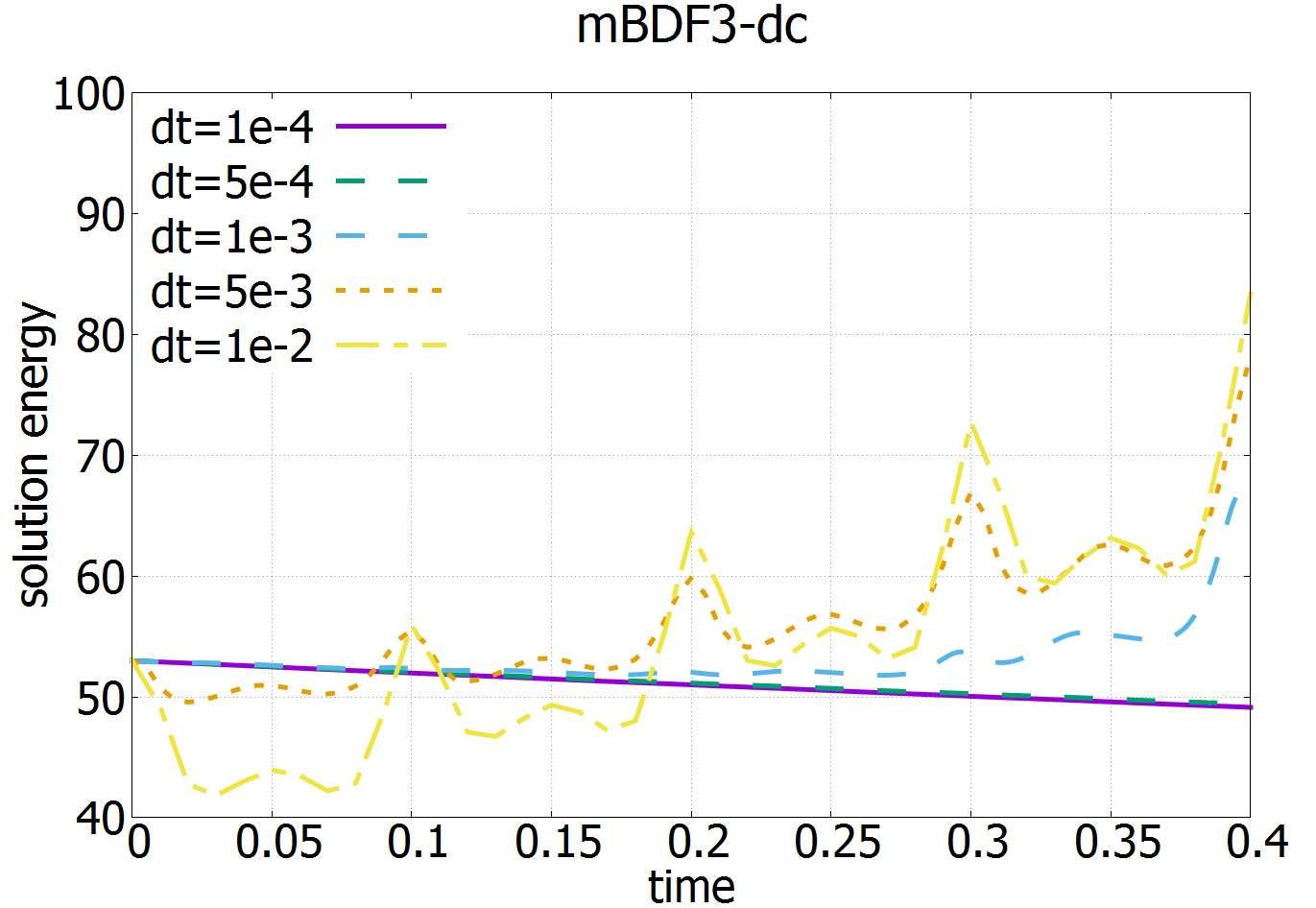}
\caption{The predicted $L^2(\Omega(t))$ norms of discrete solutions for different time steps and different schemes. Grid velocity is piecewise constant in time calculated by the first proposed approach (\ref{const-grid-vel}).}
\label{stab-0}
\end{figure}

\begin{figure}
\center \includegraphics[scale=0.14]{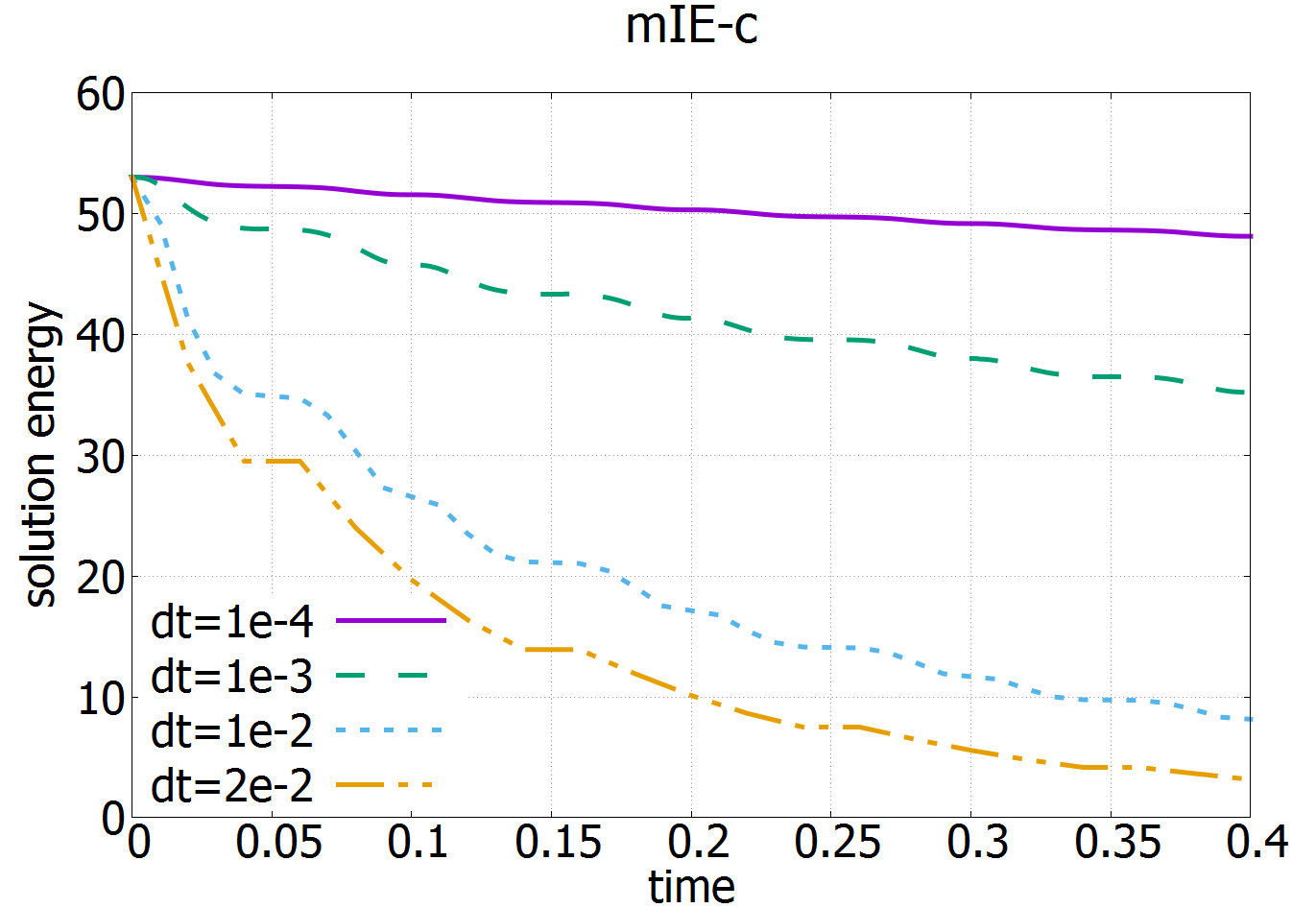}\includegraphics[scale=0.14]{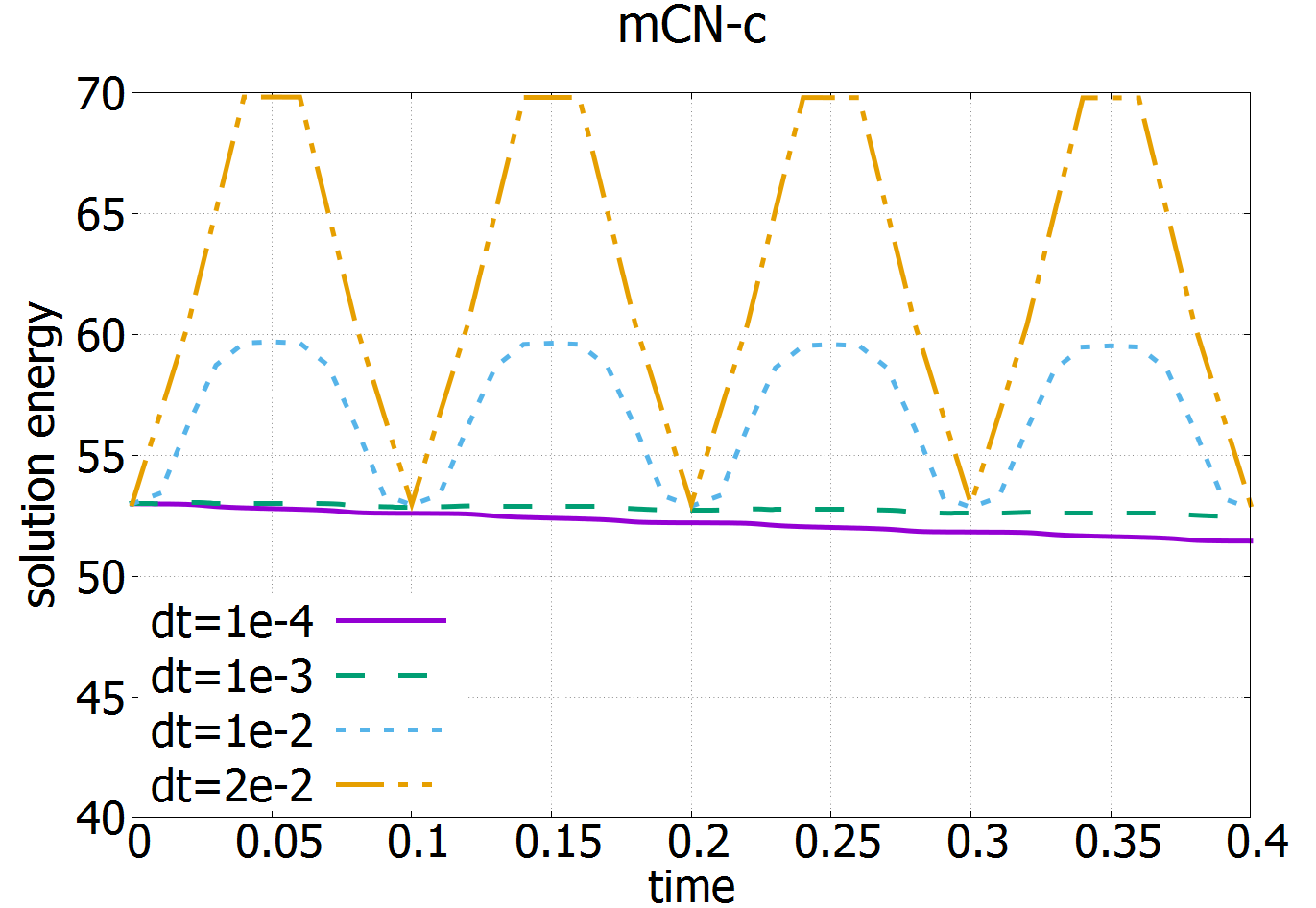}

\vspace{0.5cm}

\center \includegraphics[scale=0.14]{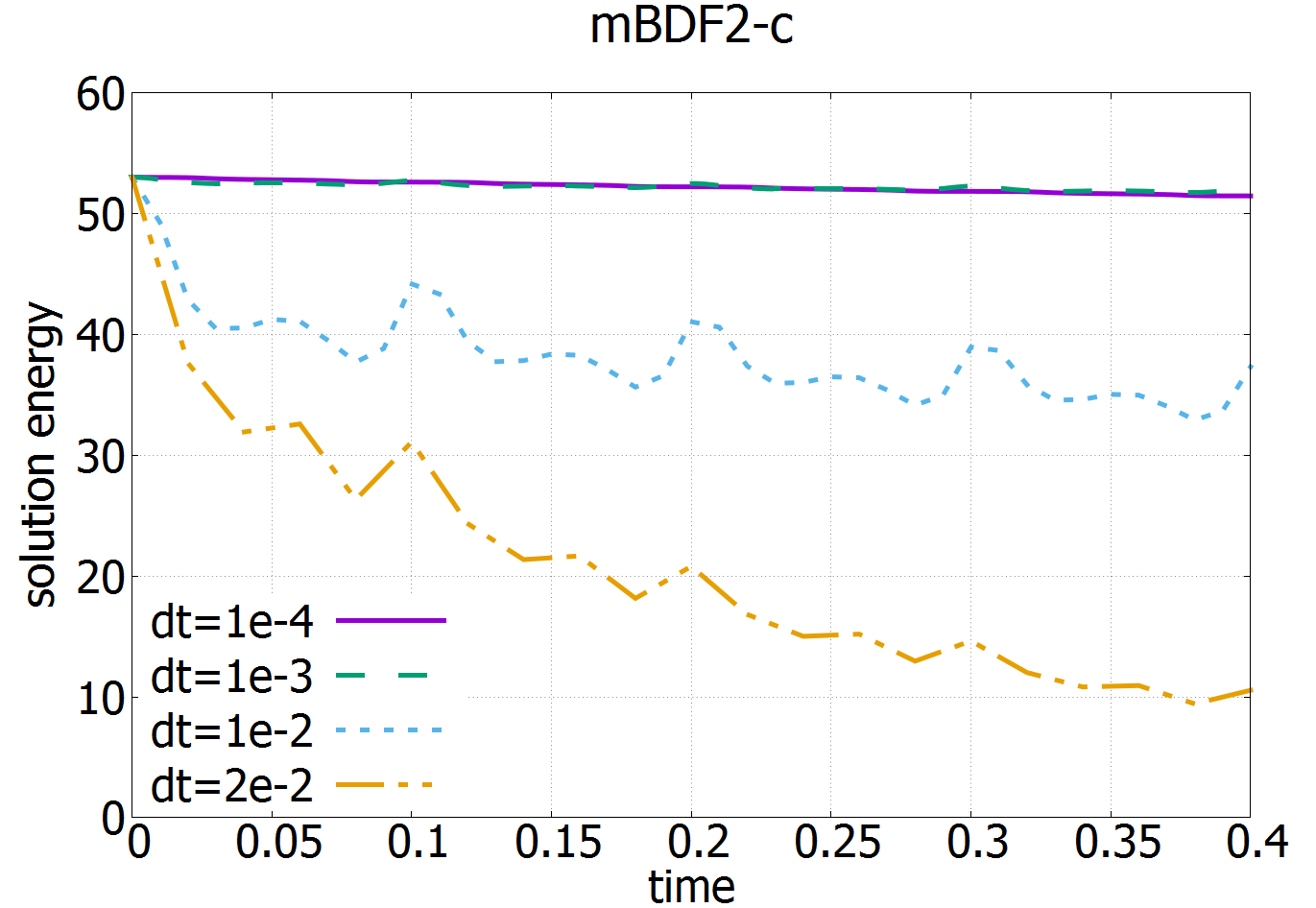}\includegraphics[scale=0.14]{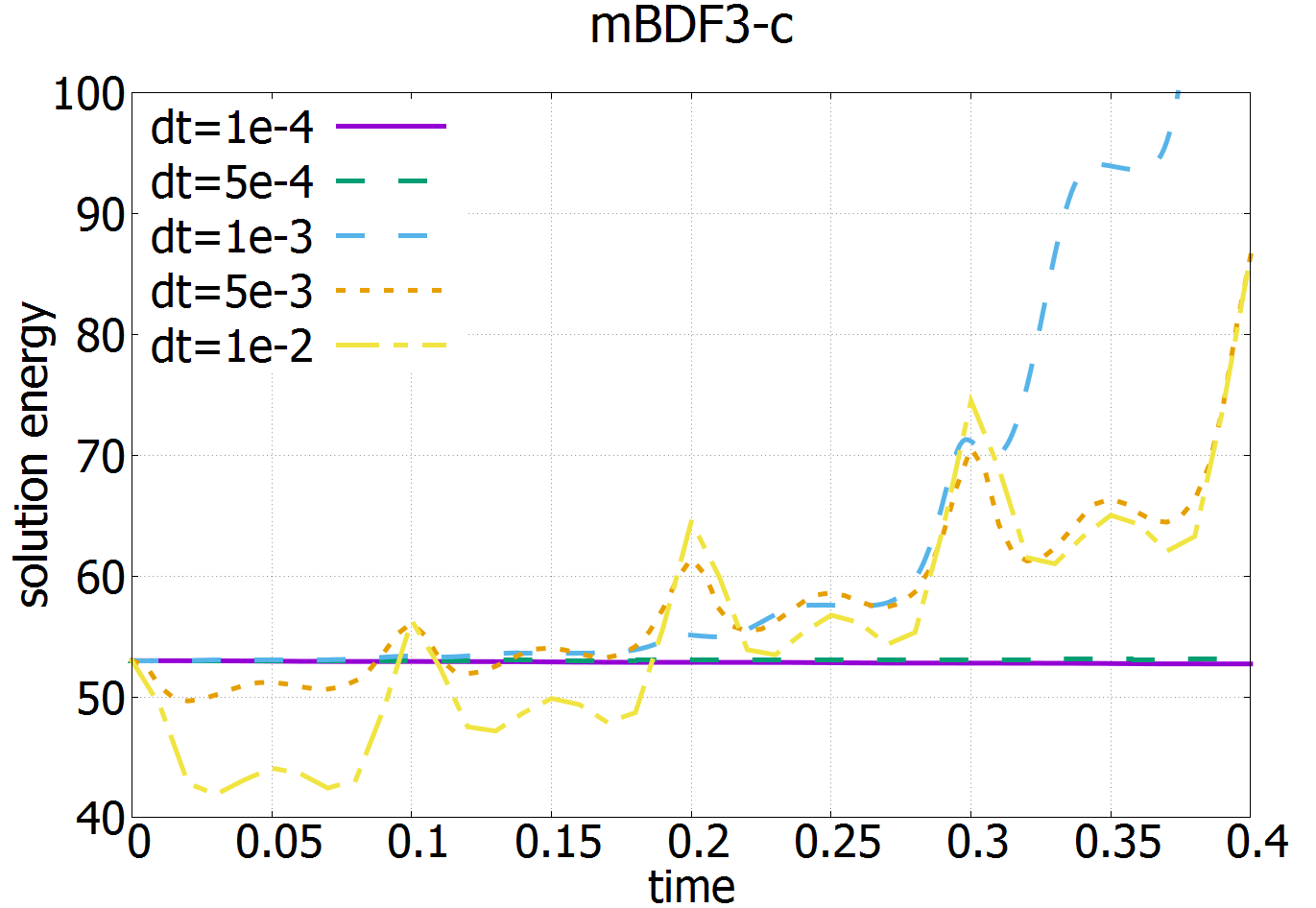}
\caption{The $L^2(\Omega(t))$ norm of discrete solution for different time steps and different schemes. Grid velocity is continuous in time calculated by the secondly proposed approach (\ref{con-grid-vel}).}
\label{stab-1}
\end{figure}

We have used the same arrangement for the mesh density and time step lengths as in the original papers \cite{formaggia-nobile04,boffi-gastaldi}. Grid velocity is calculated according to the schemes described in Section~\ref{sec:vanishSCL}. The results for stability are shown in Figures \ref{stab-0} and \ref{stab-1} for the implicit Euler scheme (mIE, m denoting {\it} modified), Crank--Nicolson (mCN), BDF2 (mBDF2) and BDF3 (mBDF3) schemes.

The firstly proposed scheme with piecewise constant in time grid velocity coincides with the method for the velocity calculation in \cite{boffi-gastaldi} and the numerical results are in exquisite arrangement with theirs. It can be noticed that, if time step is chosen sufficiently small, all schemes produce solutions with the decreasing norms as expected in the continuous case. For the cases with (relatively) large time steps, only the implicit Euler scheme preserves the decreasing behavior of the norm of solution. 

In the second approach for the grid velocity calculation, results follow the same pattern, however, the difference is noticeable for higher order methods. In this case, displacement is quadratic in time, and tests were performed in \cite{formaggia-nobile04}. Their results seem close to the ones obtained for the piecewise linear displacement, while for ours the difference seems more noticeable. Figures \ref{stab-0} and \ref{stab-1} should be compared. One can notice that for the small time steps this second approach results a smaller rate of drop of the solution energy (apart from the Euler scheme). The reason for that is still unclear, but it seems that there is less numerical diffusion in comparison with the first approach.

The second order schemes give rise to some wiggles for large time steps, while for smaller time steps the solution behaves as theoretically predicted. In \cite{boffi-gastaldi} the oscillations for Crank--Nicolson scheme are predicted on theoretical background. Thus, once again it is confirmed that application of a scheme not violating the SCL alone is not sufficient to retain stability, as noticed by various authors.

BDF3 seems more unstable for large time steps than the other considered schemes. For the cases with small time steps, the scheme stabilizes and the results are in agreement with the expected behavior. The reason for the gained instabilities most probably lies in the relation between the time step and the grid velocity, but the relation is not very clear and further investigation needs to be done. 

\subsection{Convergence}

For convergence analysis we again consider the benchmark problem posed in \cite{formaggia-nobile04}:
\begin{equation}
\begin{split}
\partial_t u - 0.1\Delta u & = f \textrm{ in }\Omega(t)\times(0,T)\\
u & = 0\textrm{ on }\partial\Omega(t)\times(0,T)\\
u(0) & = 16x(1-x)y(1-y)\textrm{ in }\Omega(0)
\end{split}
\end{equation}
with $\Omega(0)=[0,1]^2$ and the prescribed ALE map 
\begin{equation}
\kALE_t(\kbfx) = (2-\cos 10\pi t)\kbfx\textrm{ in }\kOmega=\Omega(0).
\end{equation}
The forcing term $f$ has been chosen so that the corresponding exact solution is
\begin{equation}
\ku(\kbfx,t)= 16\left( 1+\frac{1}{2}\sin(5\pi t) \right)\kx(1-\kx)\ky(1-\ky).
\end{equation}
Problem is discretized with $\mathbb{P}_2$ elements.

We have taken a sequence of decreasing time steps $0.05$, $0.01$, $0.005$, $0.001$ and computed the $L^2$--norm of the error at time $t=0.3$ over the physical domain $\Omega(0.3)$ and plotted the error against the discrete time step taken in a log--log scale.

\begin{figure}
\center 

\includegraphics[scale=0.14]{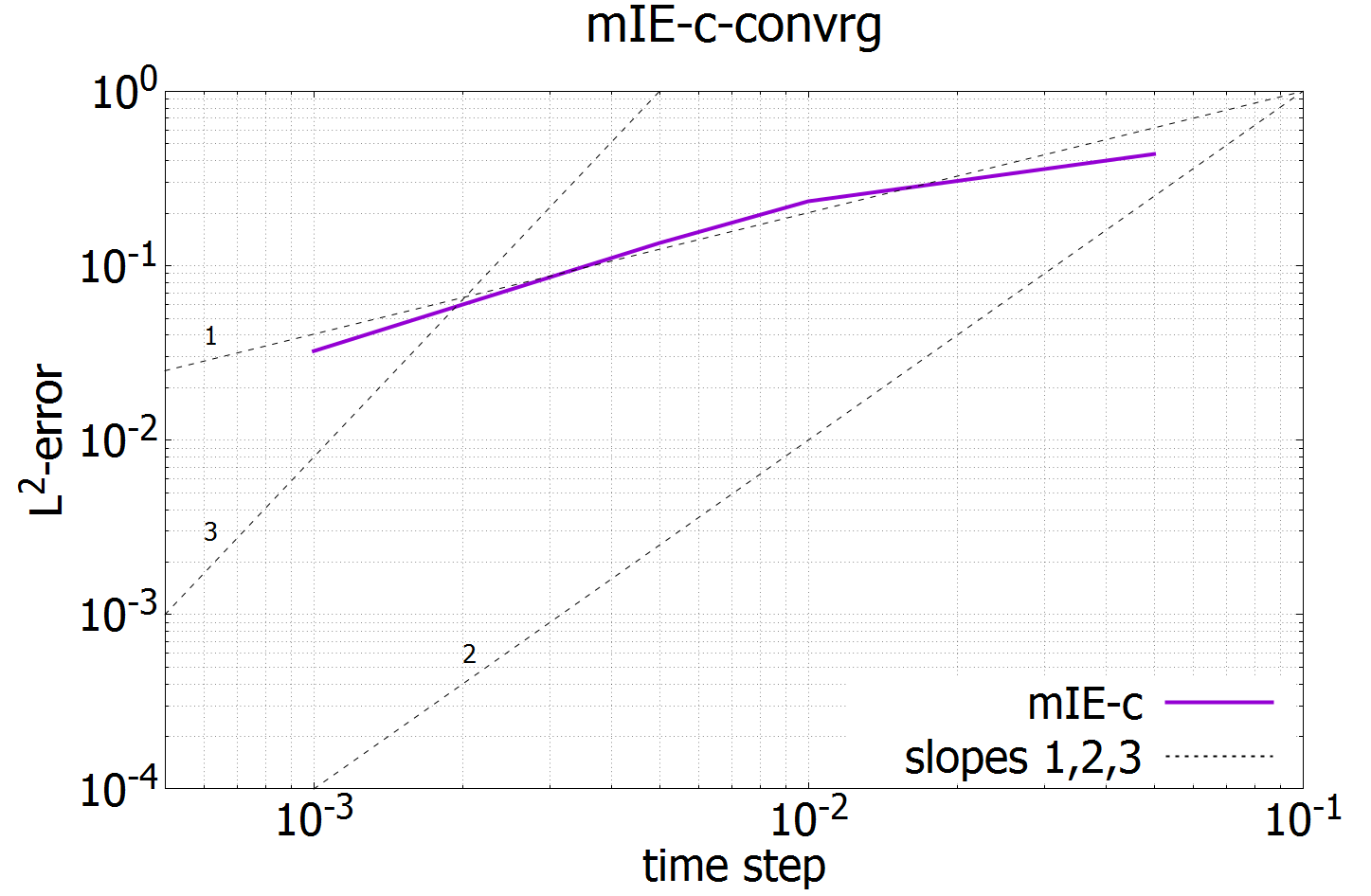}\includegraphics[scale=0.14]{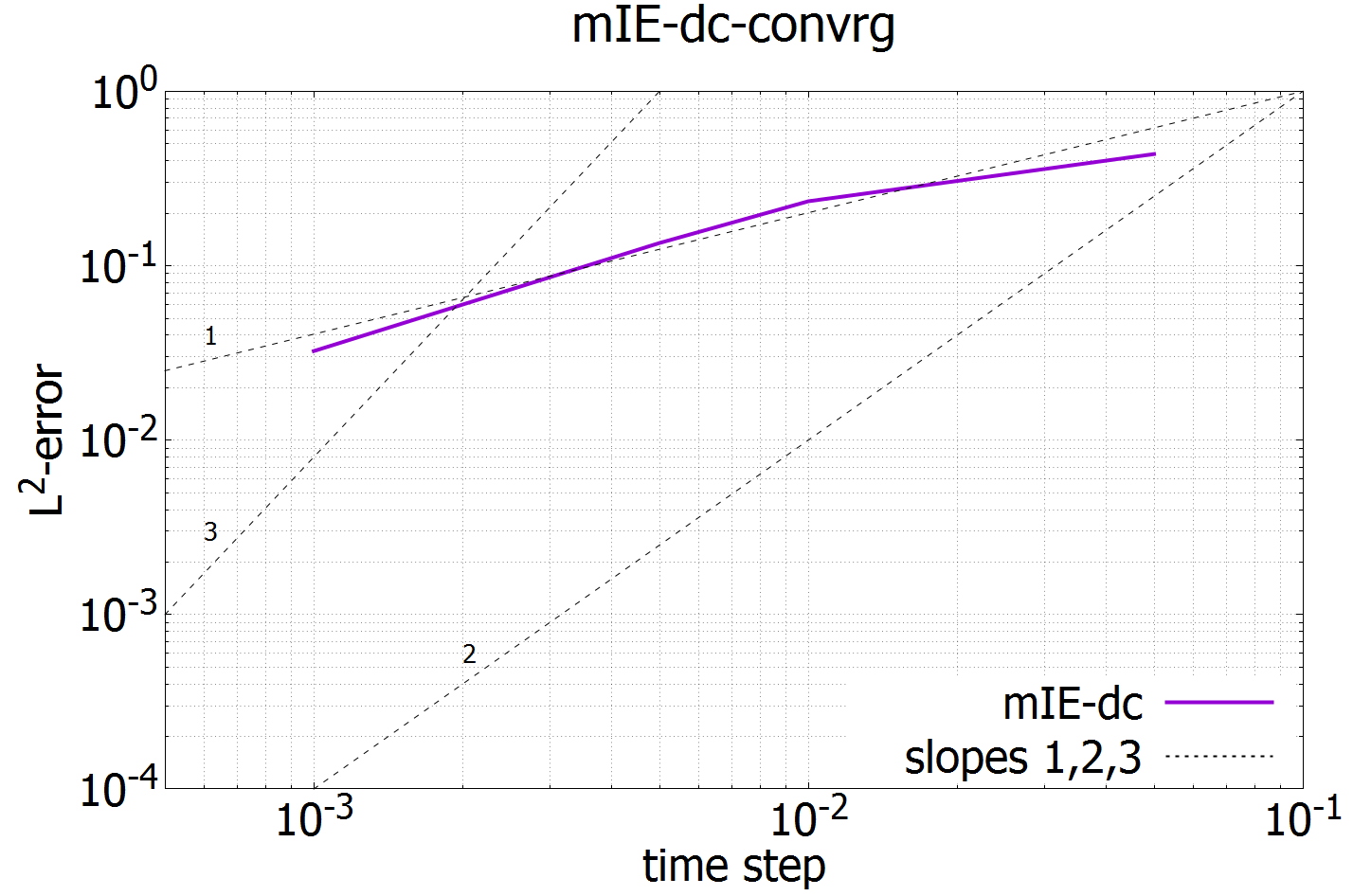}

\vspace{0.5cm}

\includegraphics[scale=0.14]{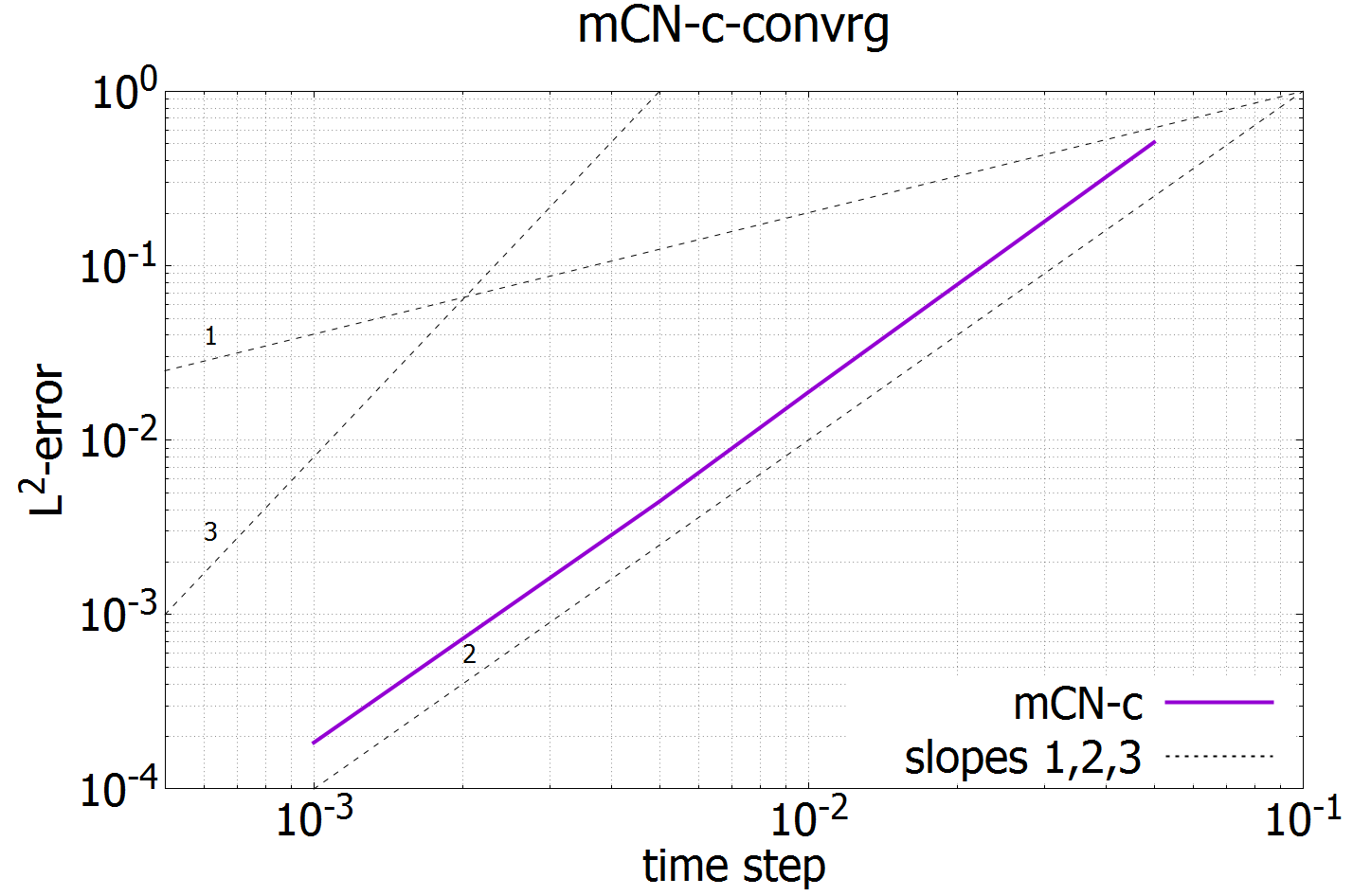}\includegraphics[scale=0.14]{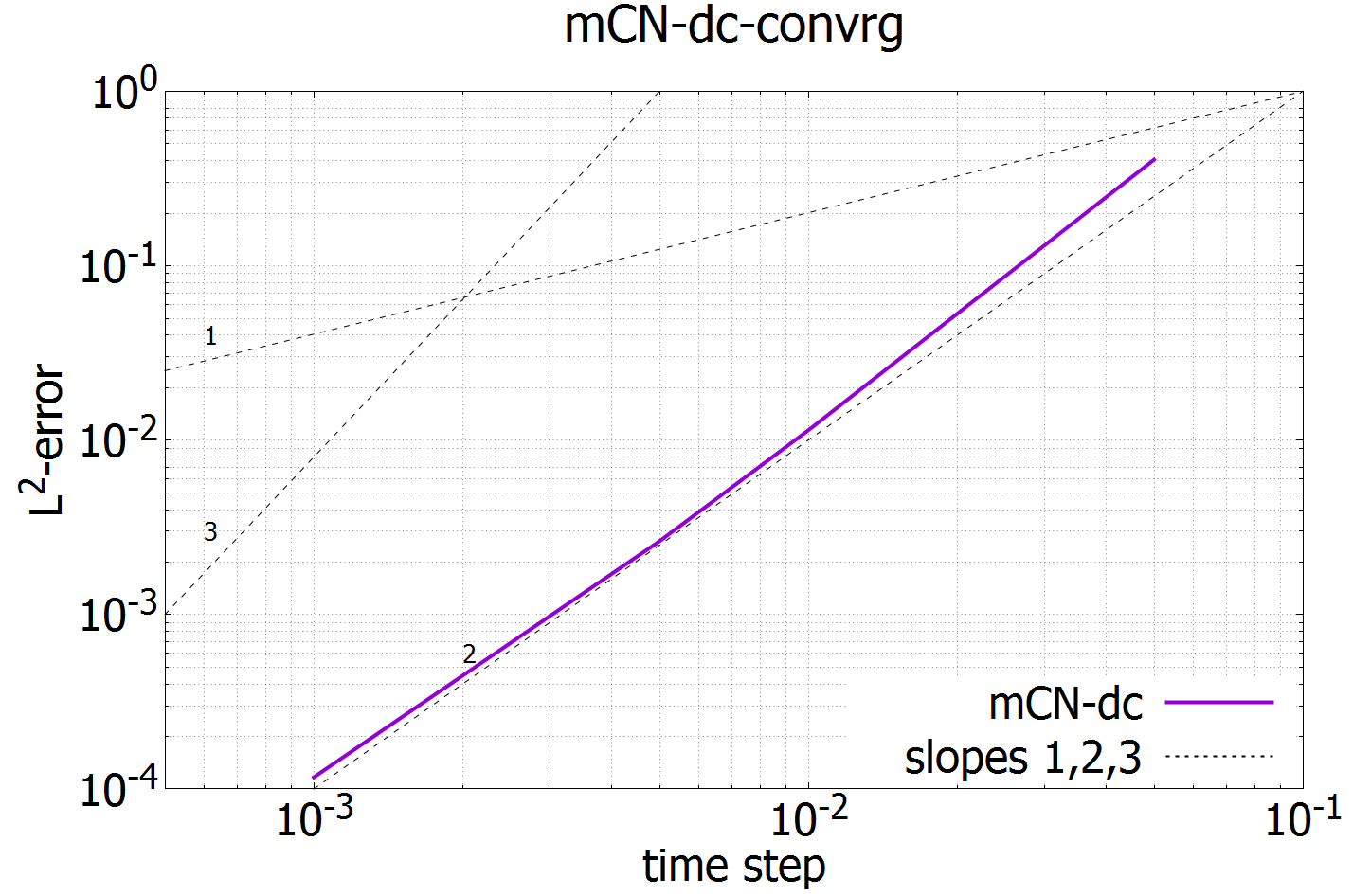}

\vspace{0.5cm}

\includegraphics[scale=0.14]{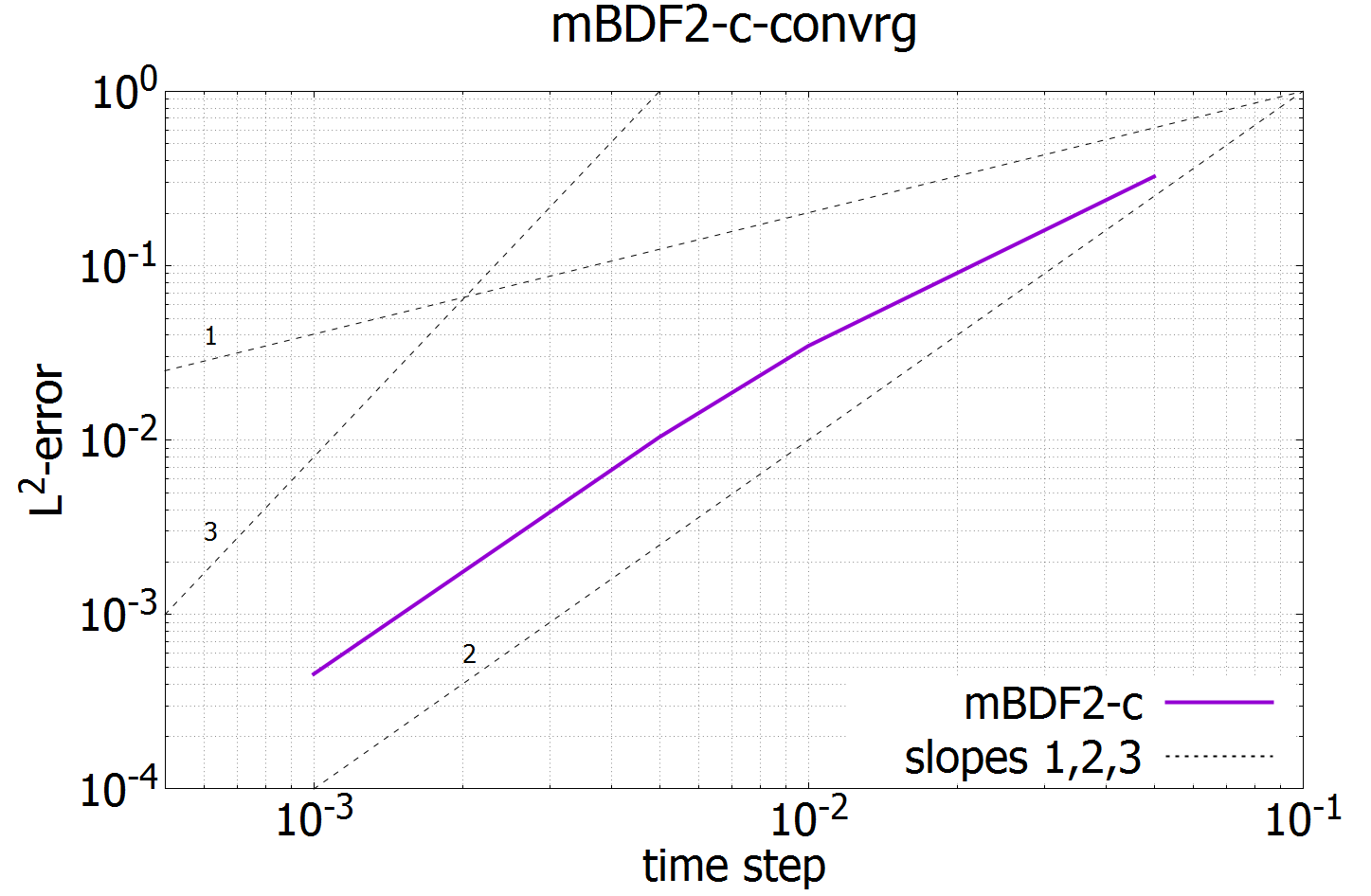}\includegraphics[scale=0.14]{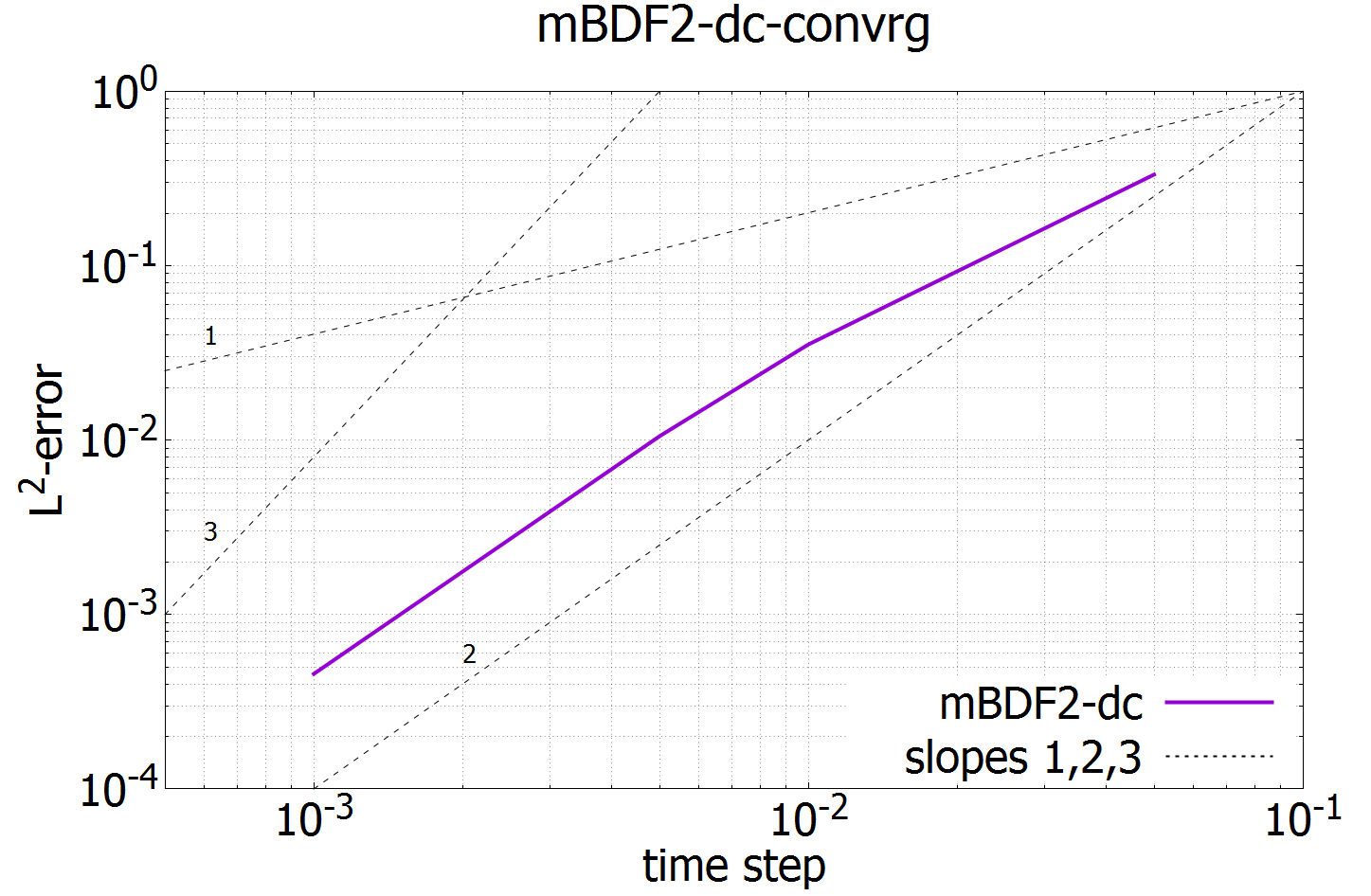}

\vspace{0.5cm}

\includegraphics[scale=0.14]{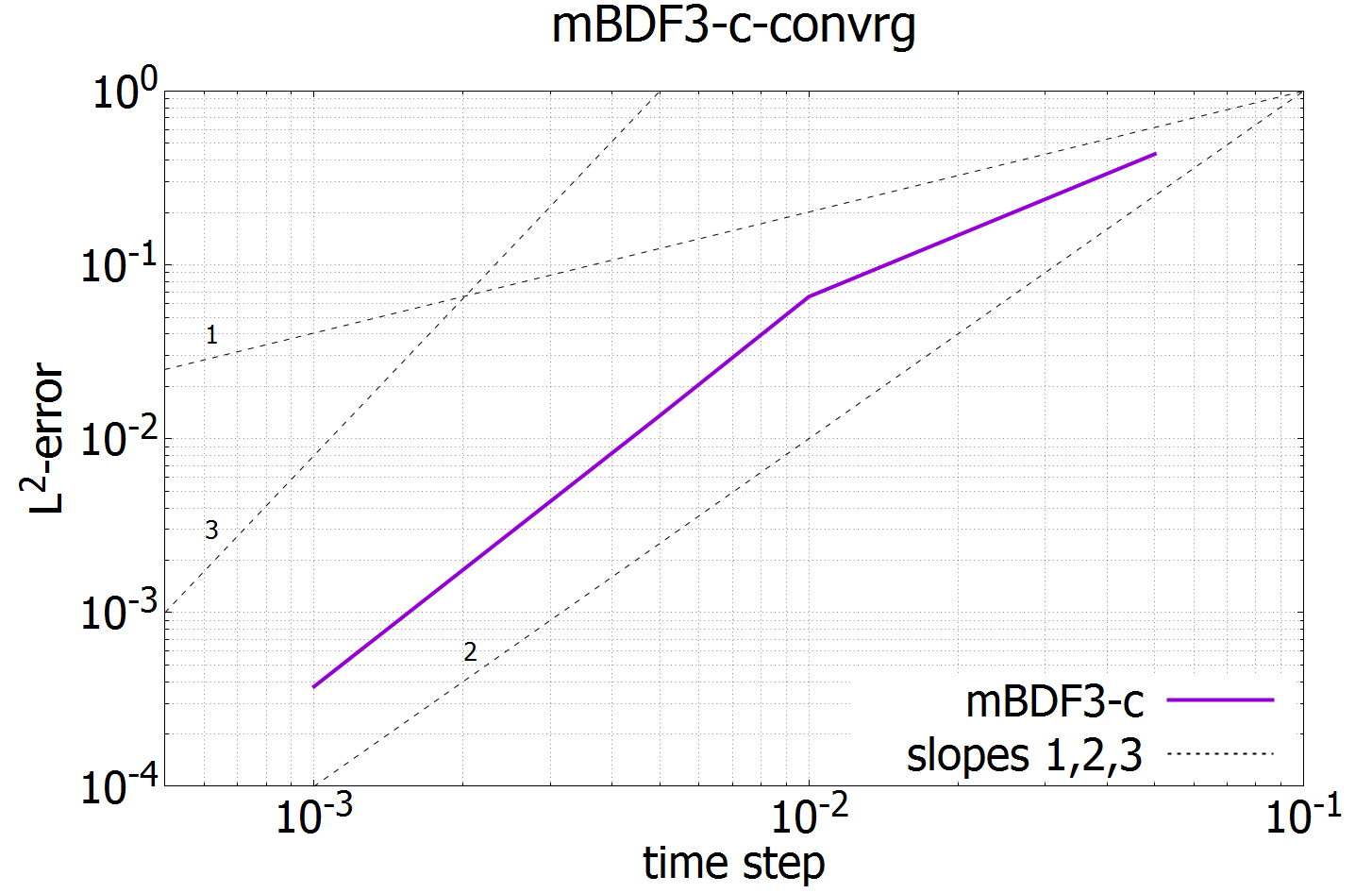}\includegraphics[scale=0.14]{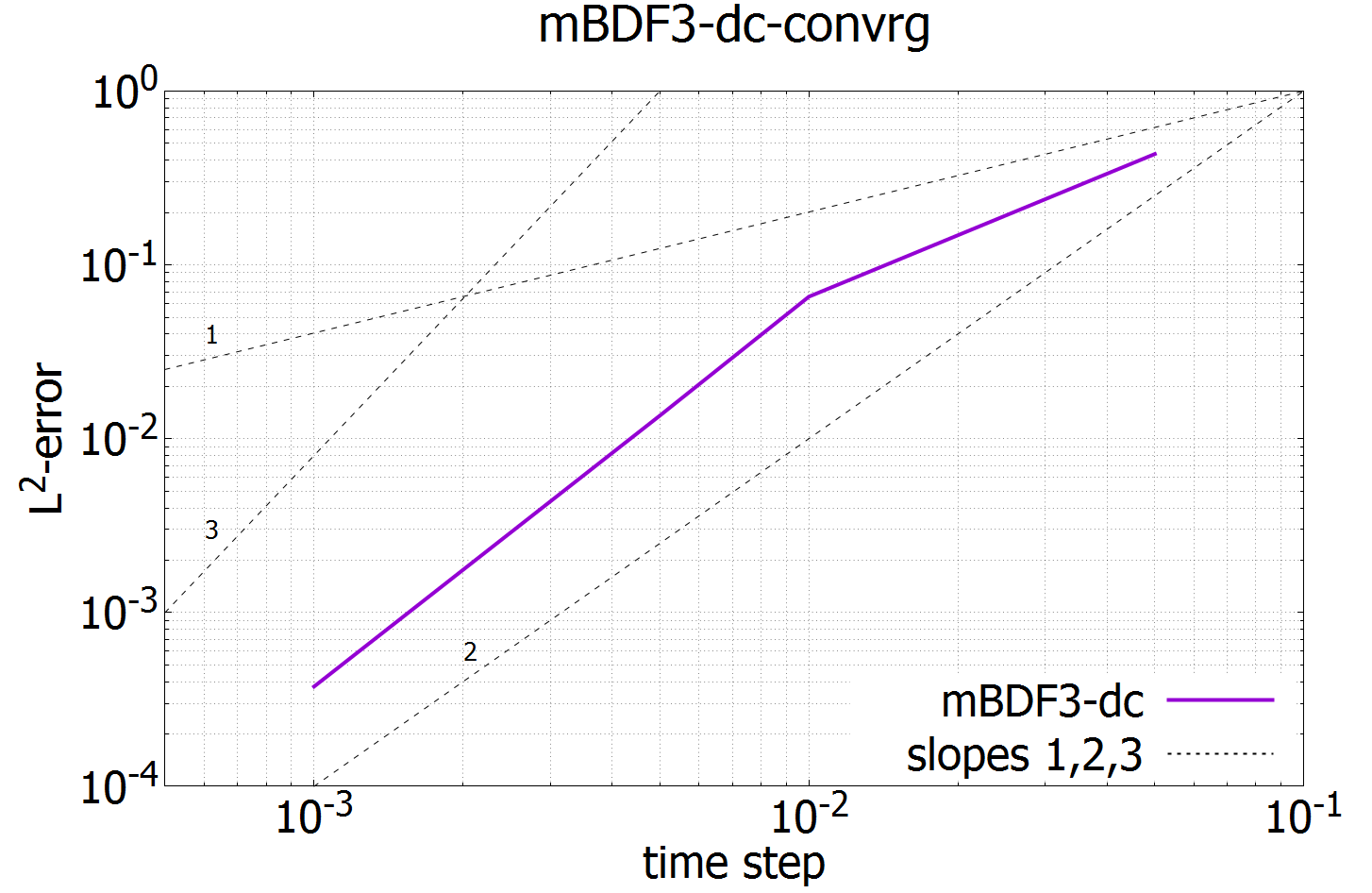}

\caption{Rates of convergence for different time--stepping schemes in log--log scale: $x$--axis represents the time step $\Delta t$ in discretization scheme ($\Delta t=0.001,0.005,0.01,0.05$), $y$--axis represents the $\Vert u^{n+1}_h-u_h(t_{n+1}) \Vert_{L^2(\Omega(t_{n+1}))}$ with $n+1$ such that $t_{n+1}=0.3$. Dashed black line denotes the slope.}

\label{rateOfCOnvrg}
\end{figure}

We observe that (see Figure \ref{rateOfCOnvrg}), apart from the BDF3, all different schemes preserve the expected order of convergence for both selections of computations of the grid velocity. Euler's scheme remains linearly convergent, while Crank--Nicolson and BDF2 schemes remain quadratically convergent. BDF3's rate of convergence is, however, between two and three, but more closer to two even for denser grids and higher order polynomials (we tested the case for spatial discretization with $\mathbb{P}_3$ elements).

\subsection{Accuracy}

To show that the  expected accuracy is preserved on moving grids, we consider a heat equation
\begin{equation}
\partial_t u-\alpha\Delta u = f\textrm{ in }\Omega(t)\times(0,T)
\end{equation}
with $\alpha=0.1$, $T=2$, $\Omega(t)=\Omega_0=[0,1]^2$, $\forall t\in(0,T)$, i.e. the domain is fixed in time. In order to simulate the moving domain problem, we prescribe the ALE map which deforms the grid (with $\kOmega=\Omega(0)$), but the domain boundary is kept unchanged. This approach allows us to compare the results obtained on the fixed grid (for which the accuracy of the schemes is known) with those obtained on moving grids.  The source term $f$ and the initial condition $u(0)$ in the above equation are chosen such that \[ u(\bfx,t)=\sin t\cos(2(x-\frac{1}{2})^2+2(y-\frac{1}{2})^2) \]
is an exact solution. Clearly, due to grid movement, the spatial discretization changes and possibly influences the accuracy of numerical solution. On dense grids, difference in accuracy due to spatial discretization should be minimally affected. As shown in the previous Section~\ref{sec:stability}, stability might play a part as well. If the produced errors for fixed and moving grid cases exhibit the same pattern (with possibly small difference due to different grids), then the accuracy is the same for both methods.

\begin{figure}
\center
\includegraphics[scale=0.28]{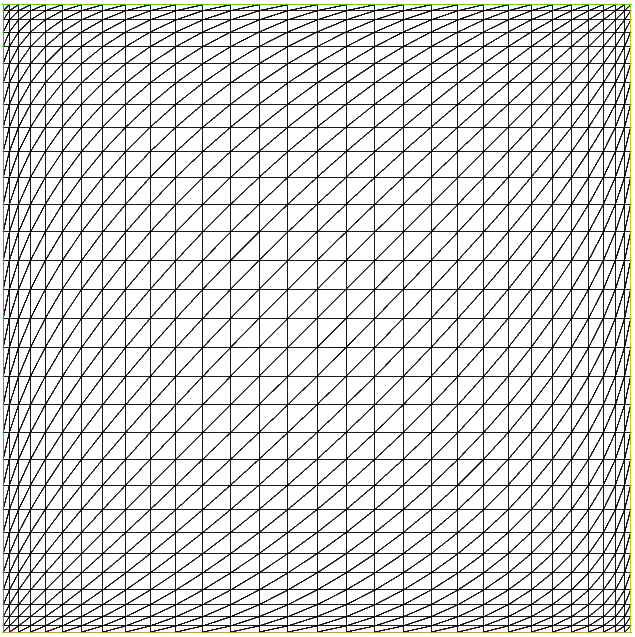}\hspace{0.3cm}\includegraphics[scale=0.28]{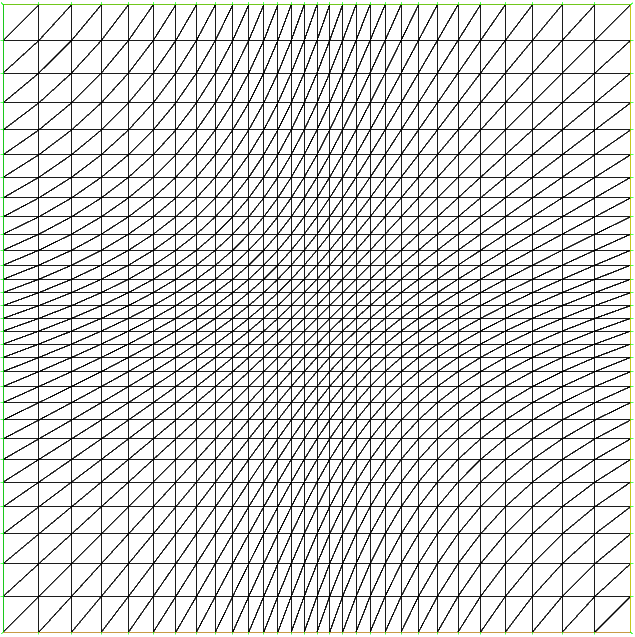}

\vspace{0.5cm}

\includegraphics[scale=0.28]{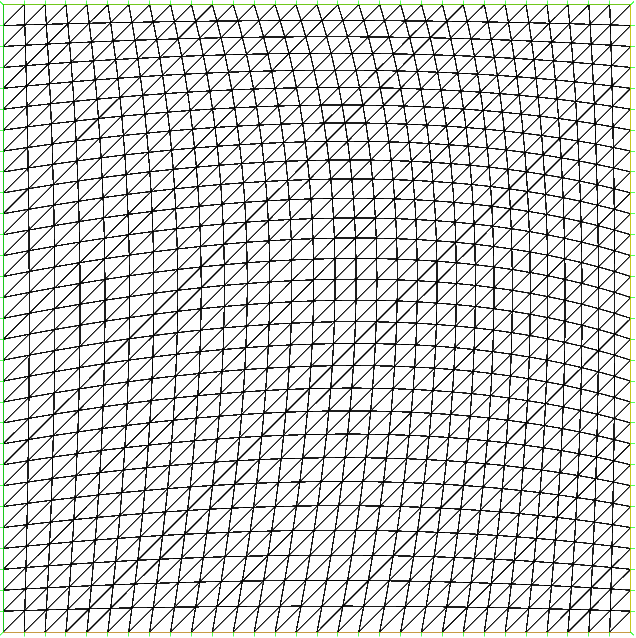}\hspace{0.3cm}\includegraphics[scale=0.28]{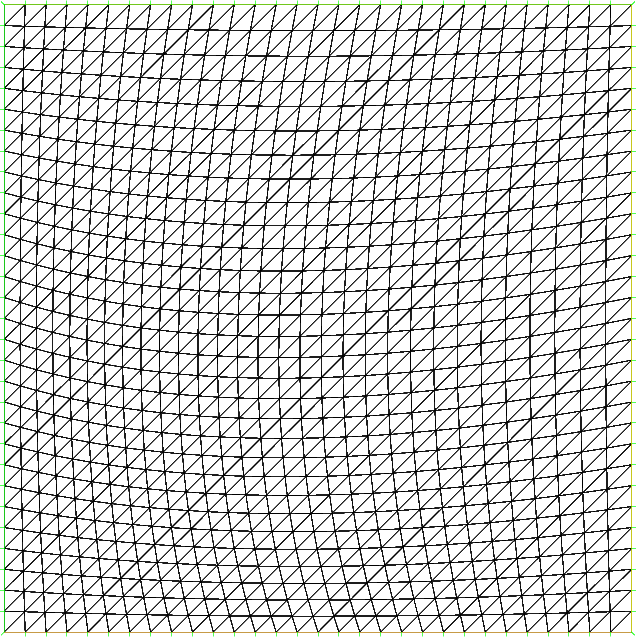}
\caption{Illustration of the action of the ALE map ($\kALE^A$ top, $\kALE^B$ bottom). The initial grid is uniform on the domain $[0,1]^2$.}
\label{grid-ill}
\end{figure}

\begin{figure}

\center 

\includegraphics[scale=0.14]{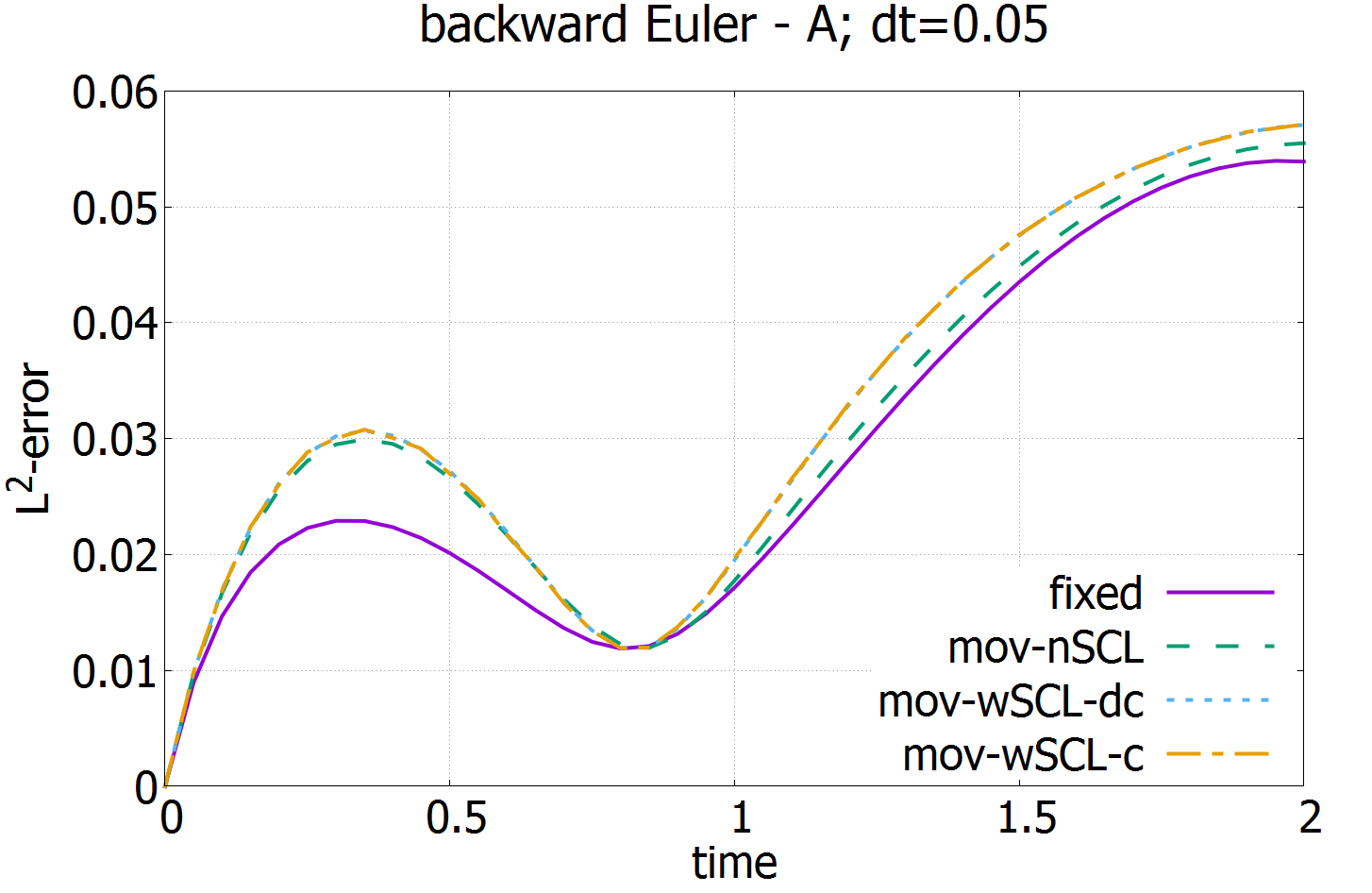}\includegraphics[scale=0.14]{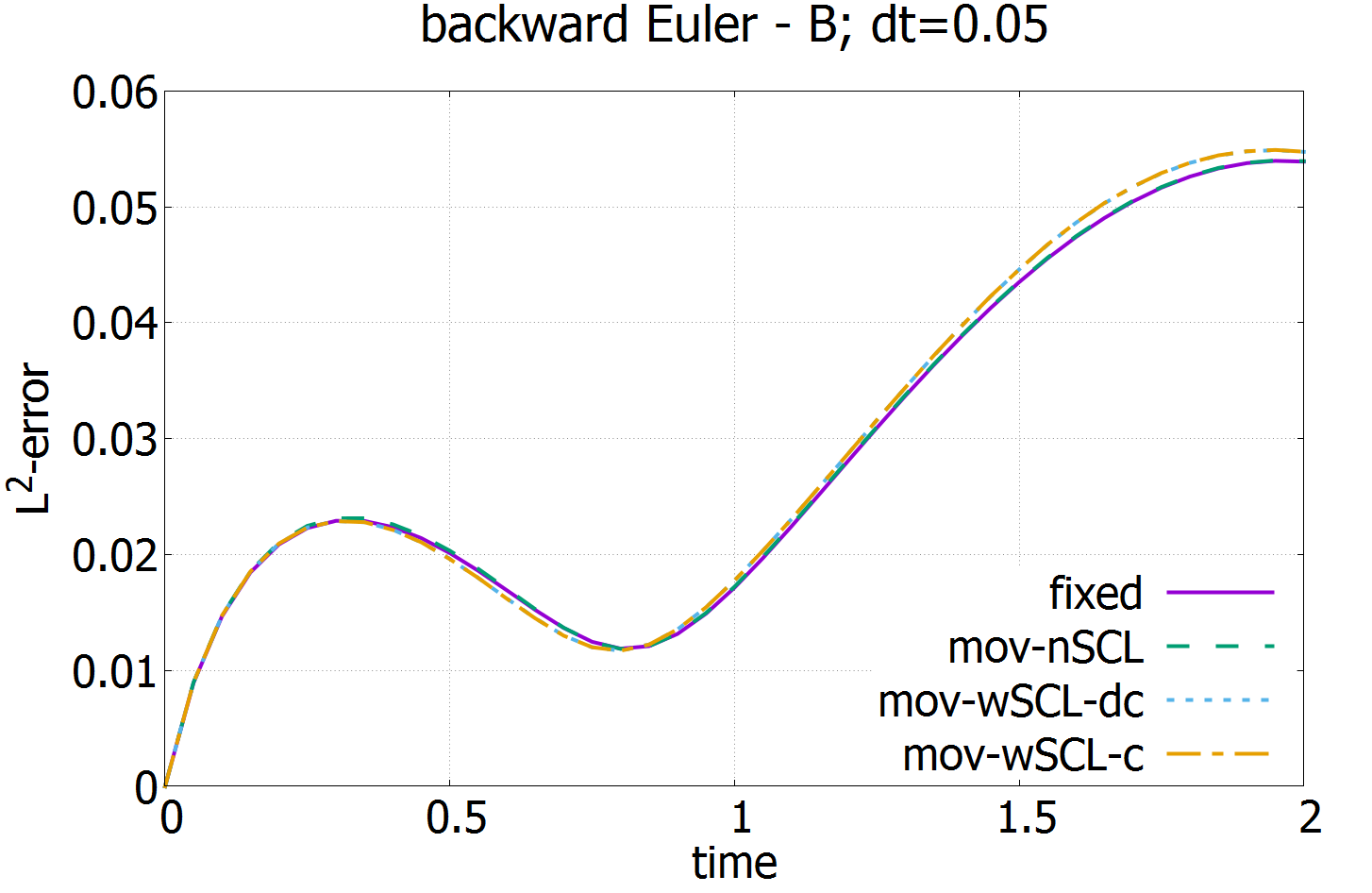}

\vspace{0.4cm}

\includegraphics[scale=0.14]{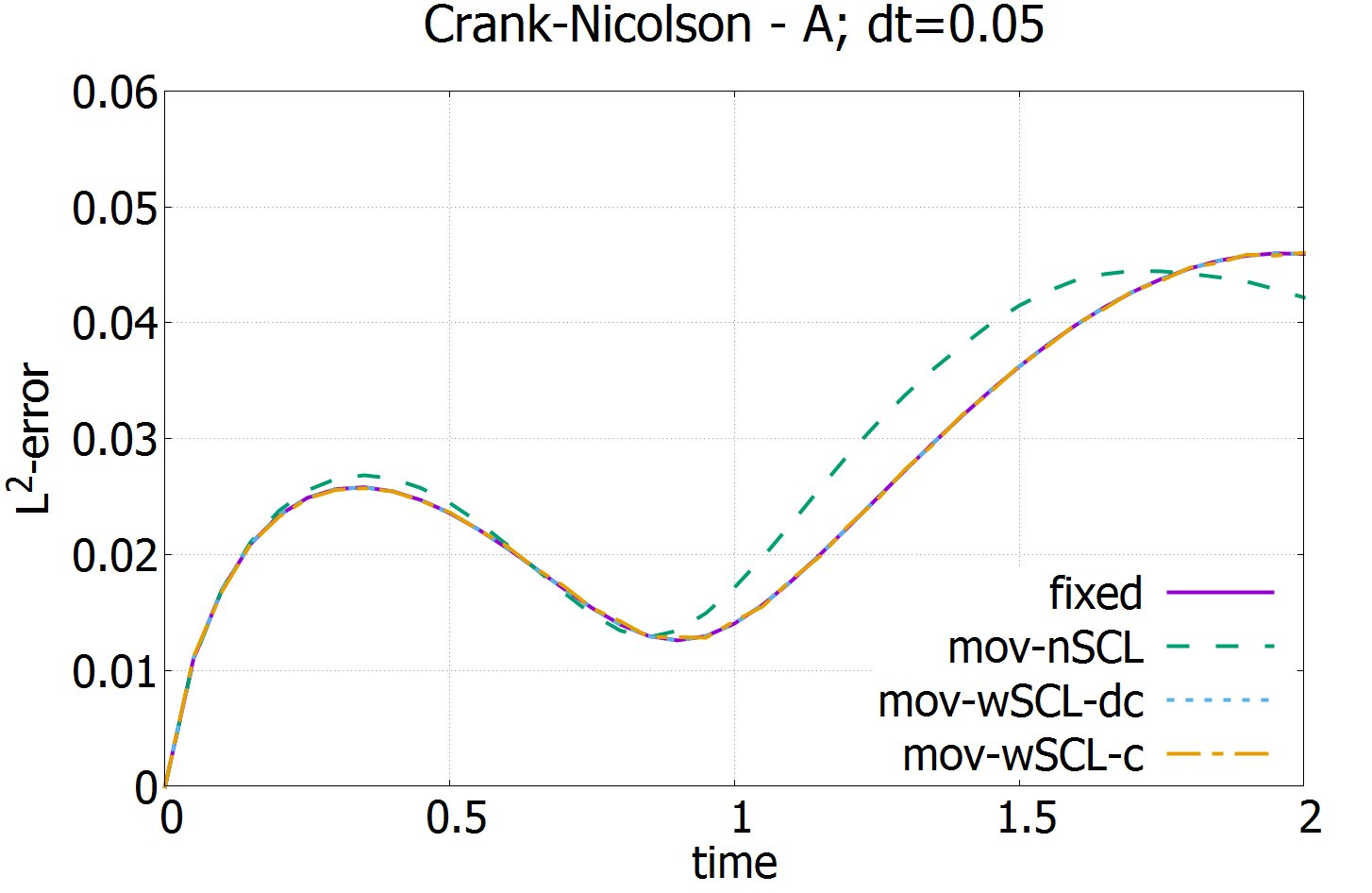}\includegraphics[scale=0.14]{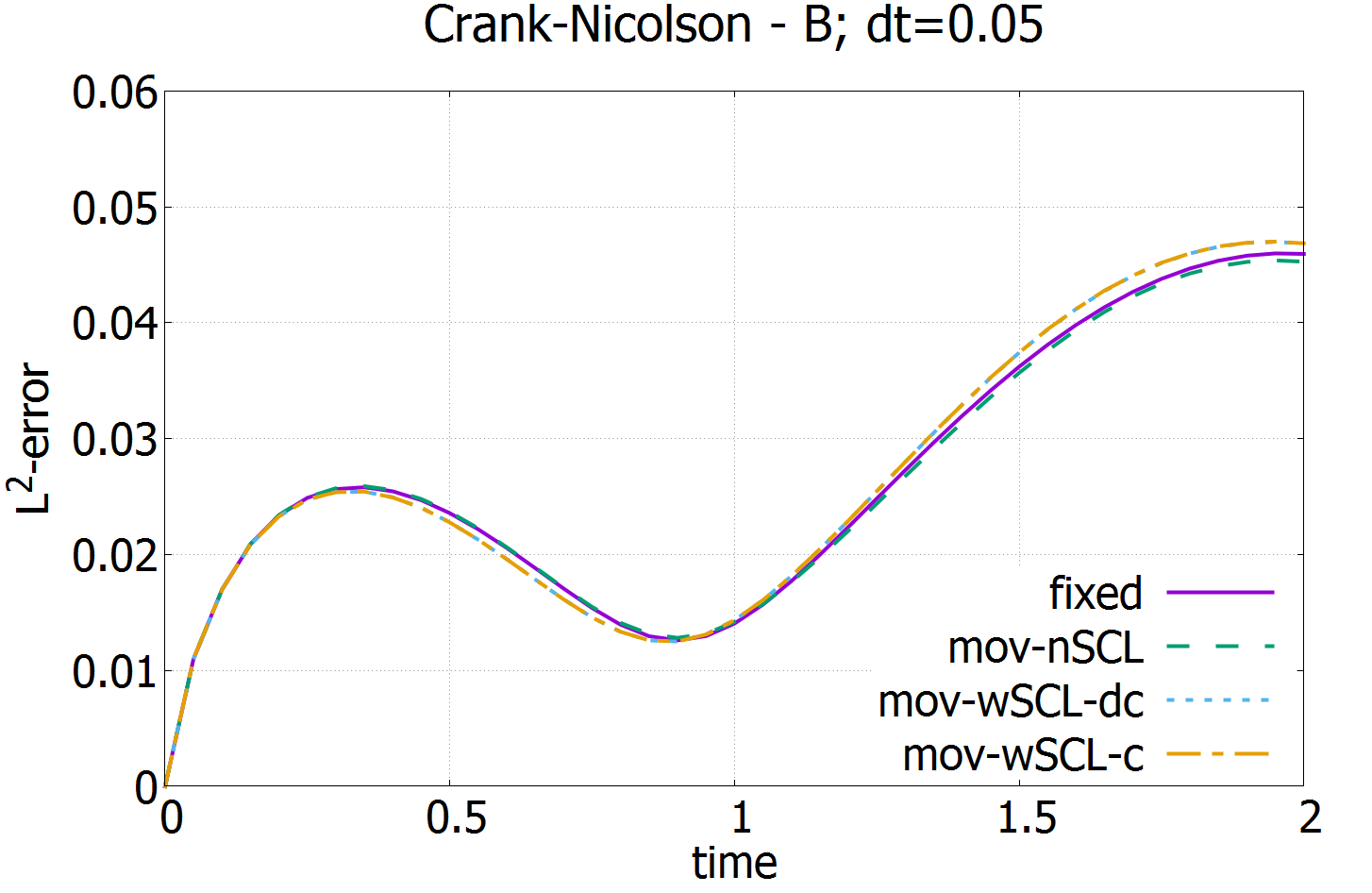}

\vspace{0.4cm}

\includegraphics[scale=0.14]{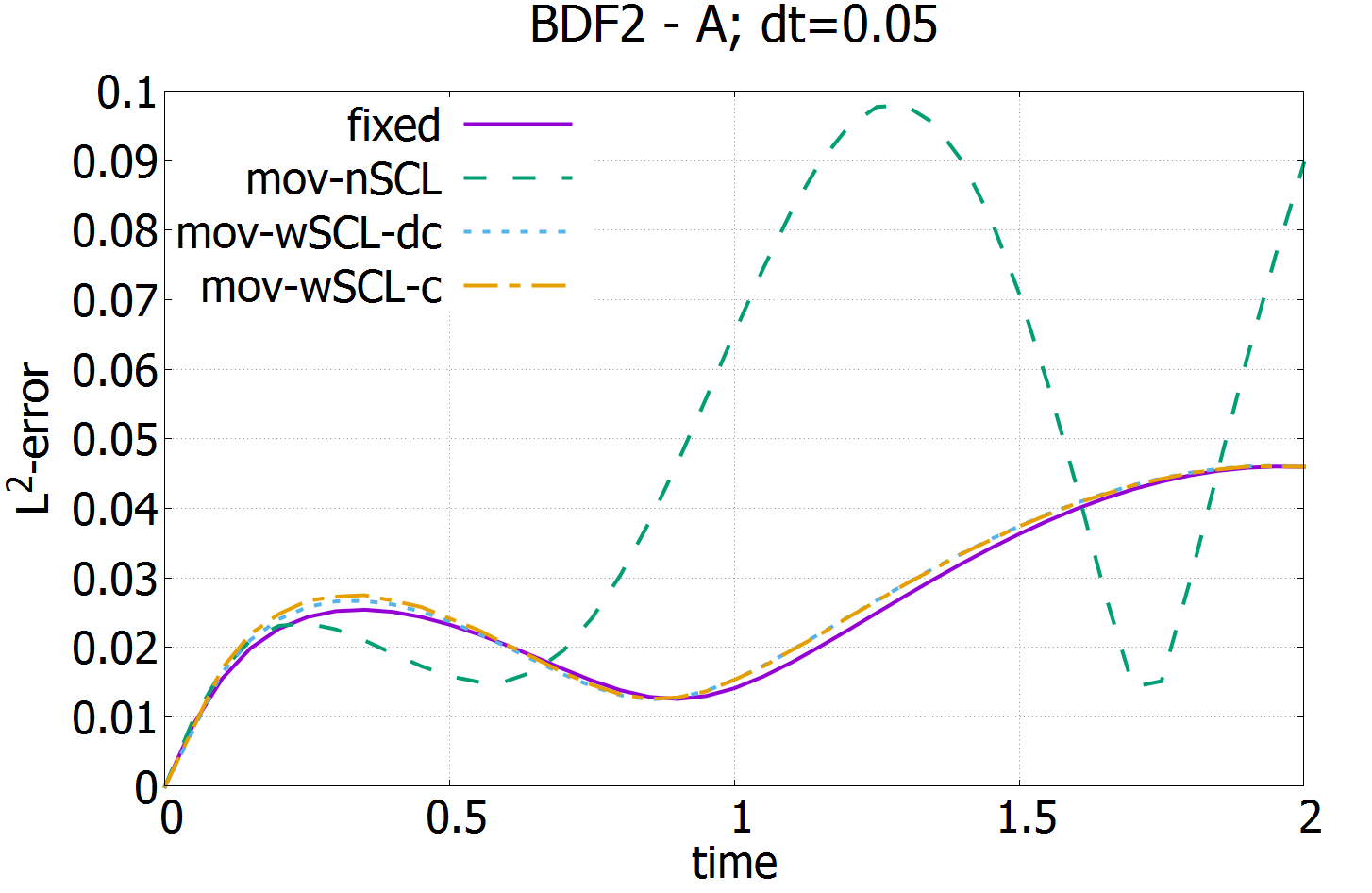}\includegraphics[scale=0.14]{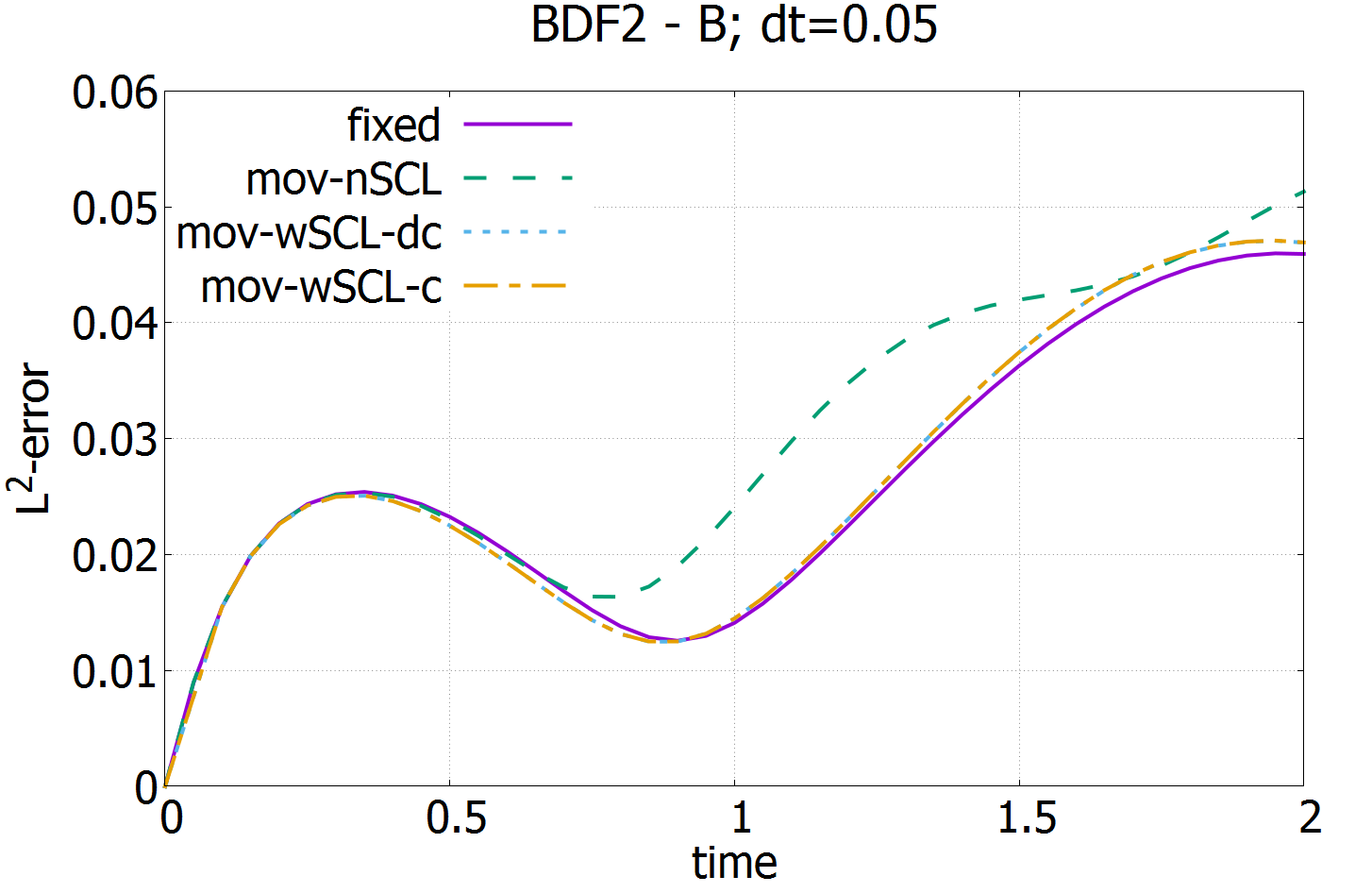}

\vspace{0.4cm}

\includegraphics[scale=0.14]{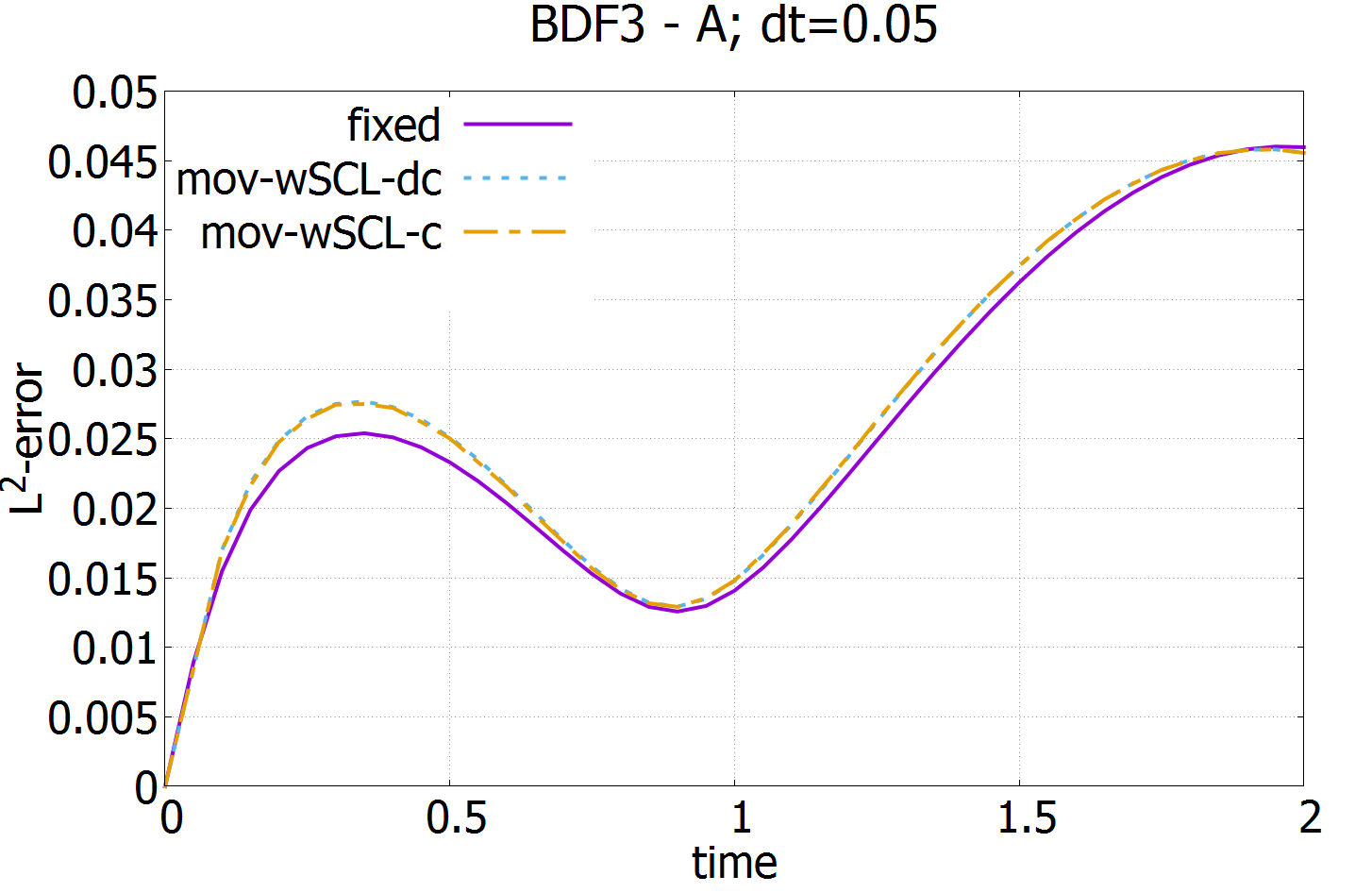}\includegraphics[scale=0.14]{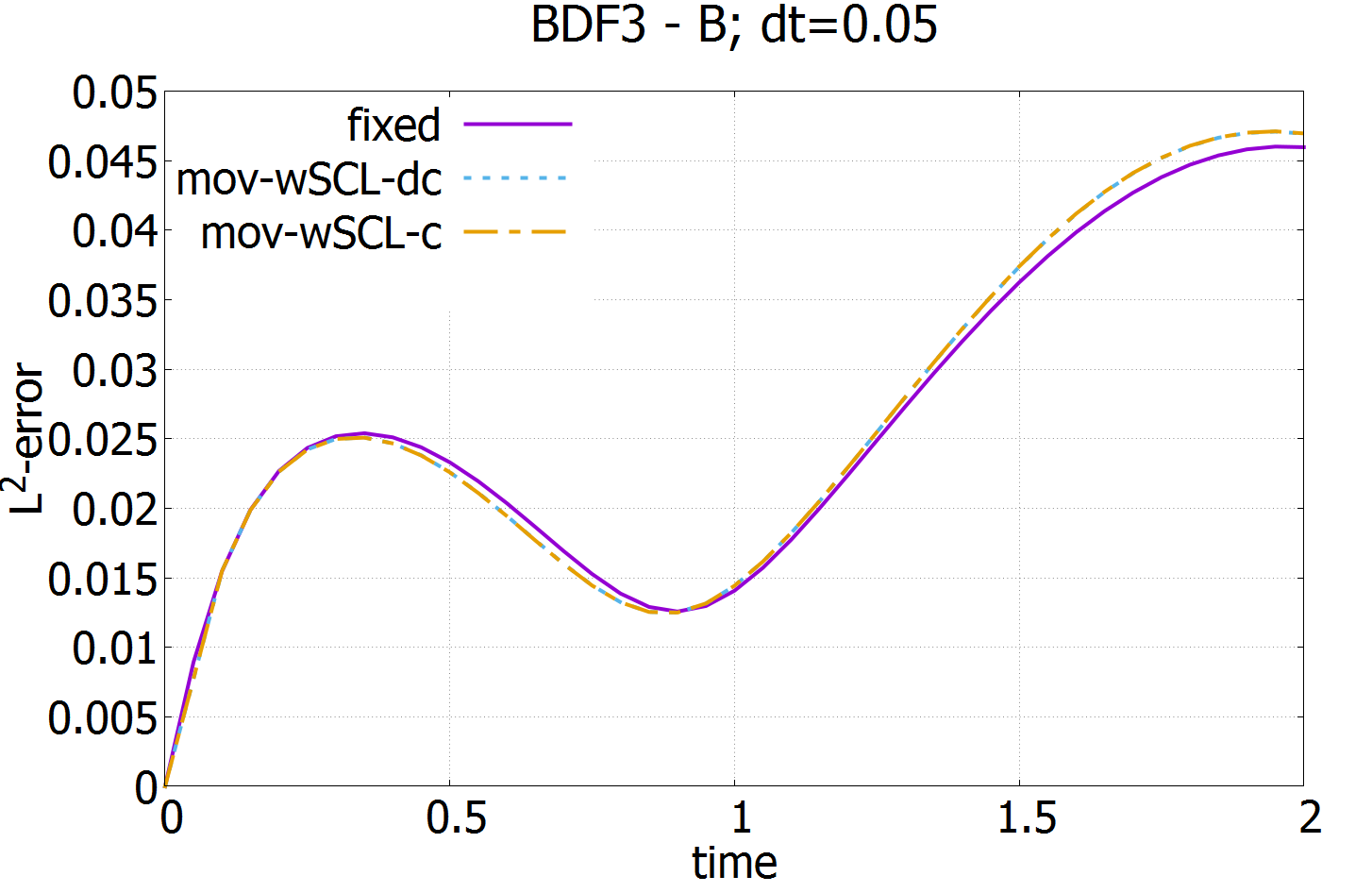}
\caption{The $L^2$ norms of errors between the exact and numerical solutions for different schemes and time step $\Delta t=0.05$. In the legend "{\it fixed}" refers to solutions obtained on fixed grids, while "{\it mov-wSCL}" and "{\it mov-nSCL}" for ones obtained on moving grids with the proposed non--violating SCL scheme ("$dc$" and "$c$" standing for the first and second approach in grid velocity calculation), and the classical, SCL--violating scheme, respectively. }
\label{fig:acc}
\end{figure}

The grid is moved according to the ALE maps given bellow (interpolated onto $\mathbb{P}_1$ space):
\[ \kALE^A(\kbfx,t)= \begin{bmatrix}
\kx + \frac{1}{2}\sin(\pi t)\sin(\pi\kx(1-\kx)(\kx-\frac{1}{2}))\\
\ky + \frac{1}{2}\sin(\pi t)\sin(\pi\ky(1-\ky)(\ky-\frac{1}{2}))
\end{bmatrix} \]
and
\[ \kALE^B(\kbfx,t)= \begin{bmatrix}
\kx + \sin(\pi t)\kx(1-\kx)\ky(1-\ky)\\
\ky + \sin(\pi t)\kx(1-\kx)\ky(1-\ky)
\end{bmatrix}. \]
The above maps are constructed with the objective to change the volume of triangles in order to emphasize the problematics rising from the violation of the SCL constraint (see Figure \ref{grid-ill}).
We also plot the results obtained for SCL--violating schemes from which the newly proposed schemes were derived for an illustration. From the Figure \ref{fig:acc} we can observe that the classical SCL--violating schemes may produce an error that does not follow the pattern of the one obtained in the fixed grid.  The difference is especially noticeable for BDF2 method, although this is expected from the theoretical discussion in 5.3. The proposed SCL--non--vioalting schemes produce errors that follow the patterns of the ones obtained in the fixed grid in very good agreement. For the small time steps agreement is excellent, while for the larger time steps solutions might suffer from the artificial numerical diffusion already discussed in Section~\ref{sec:stability}.

\section{Conclusion}\label{sec:conclusion}

A modified approach of handling PDEs on time--dependent domains with finite element method  has been introduced within Arbitrary Lagrangian Eulerian framework. The approach exploits the full potential of the polynomial form of grid velocity by integrating the time integrals exactly. Consequently, the SCL identity is trivially satisfied. While much more work remains to be done on the question of stability, it seems that in case when discrete time step is sufficiently small to keep the scheme stable, the accuracy of the schemes is maintained. From the numerical results a conclusion can be made that for the problems on moving grids not violating SCL alone is not enough for stability (and convergence) of the schemes. Although it yet remains to be confirmed, it looks like the terms that do not explicitly involve the grid velocity play a more significant role than expected. We suspect this by considering the problem re--posed on the referential configuration. In that case the dependence of all terms on the ALE map (and thus on the grid motion) is more emphasized than in the case when problem is posed on current configuration. This dependence emerges in the form of Jacobians and gradients of ALE map. An attempt in more theoretical analysis will be made in the following papers.

The main advantage of the newly proposed approach is in its simplicity for the generalization to an arbitrary temporally high--order scheme without (explicitly) worrying about the discrete SCL. Independently on the chosen scheme for the discretization of temporal derivative, it is always possible to satisfy the SCL. Moreover, the satisfaction of SCL is trivial due to the construction of the formulation which is based on the differential statement of SCL. Two alternative ways of grid velocity calculation have been presented, and their influences on the scheme stability have been investigated.

A quite detailed procedure on the pullback to the referent configuration has been presented. Few of the most popular schemes have been considered and the details on their modifications within the new approach have been given.

All of the implementation has been done with the {\it FreeFEM++} software (\cite{hecht,decoene}). Numerical validation showed good agreement with the benchmark  problems introduced by other authors.

\bibliographystyle{siamplain}
\bibliography{references}

\end{document}